%% file: sarah_nonconvex_arxiv.tex
\documentclass[letterpage, 11pt, notitlepage]{article}

\usepackage[margin=1.0in]{geometry}

\usepackage{enumitem}
\usepackage{times}
\usepackage{graphicx} 
\usepackage{subfigure}

\usepackage{algorithm}
\usepackage{algorithmic}

\usepackage{hyperref}
\usepackage{breakurl}

\AtBeginShipout{%
  \ifnum\value{page}>1 %
    \typeout{* Additional boxing of page `\thepage'}%
    \setbox\AtBeginShipoutBox=\hbox{\copy\AtBeginShipoutBox}%
  \fi
}

\input{mystyle.tex}

\makeatother

\usepackage{blindtext,titlefoot}

\begin{document}

\title{\fontsize{20}{20}\selectfont Finite-Sum Smooth Optimization with SARAH}

\author{
Lam M. Nguyen
\and
Marten van Dijk
\and
Dzung T. Phan
\and
Phuong Ha Nguyen
\and
Tsui-Wei Weng
\and
Jayant R. Kalagnanam
}

\maketitle

\unmarkedfntext{Correspondence to: \textbf{Lam M. Nguyen}, IBM Research, Thomas J. Watson Research Center, Yorktown Heights, NY, USA. Email: \href{mailto:lamnguyen.mltd@ibm.com}{LamNguyen.MLTD@ibm.com}}

\begin{abstract}

The total complexity (measured as the total number of gradient computations) of a stochastic first-order optimization algorithm that finds a first-order stationary point of a finite-sum smooth nonconvex objective function $F(w)=\frac{1}{n} \sum_{i=1}^n f_i(w)$
 has been proven to be at least $\Omega(\sqrt{n}/\epsilon)$ for $n \leq \Ocal(\epsilon^{-2})$ where $\epsilon$ denotes the attained accuracy $\mathbb{E}[ \|\nabla F(\tilde{w})\|^2] \leq \epsilon$ for the outputted approximation $\tilde{w}$  \cite{fang2018spider}.
In this paper, we provide a convergence analysis for a slightly modified version of the SARAH algorithm \cite{Nguyen2017_sarah,Nguyen2017_sarahnonconvex} and achieve total complexity that matches the lower-bound worst case complexity in \cite{fang2018spider} up to a constant factor when $n \leq \Ocal(\epsilon^{-2})$ for nonconvex problems. For convex optimization, we propose SARAH++ with sublinear convergence for general convex and linear convergence for strongly convex problems; and we provide a practical version for which numerical experiments on various datasets show an improved performance.
\end{abstract}

\input{sections/introduction.tex}

\input{sections/method.tex}

\input{sections/sarahpp_experiments.tex}


\input{sections/adaptive.tex}

\input{sections/conclusion.tex}


\bibliography{reference}
\bibliographystyle{plain}

\appendix
\newpage

\input{sections/appendix.tex}

\end{document}

%% file: mystyle.tex
\usepackage{epsfig}
\usepackage{amssymb}
\usepackage{amsmath}
\usepackage{amsthm}
\usepackage{amsfonts}
\usepackage{bbding}
\usepackage{array}
\usepackage{caption,tabularx,booktabs}
\usepackage{arydshln}
\usepackage{xargs}                      
\usepackage[pdftex,dvipsnames]{xcolor}  

\usepackage[colorinlistoftodos,prependcaption,textsize=tiny]{todonotes}
\newcommandx{\unsure}[2][1=]{\todo[inline,linecolor=red,backgroundcolor=red!25,bordercolor=red,#1]{#2}}
\newcommandx{\change}[2][1=]{\todo[linecolor=blue,backgroundcolor=blue!25,bordercolor=blue,#1]{#2}}
\newcommandx{\info}[2][1=]{\todo[linecolor=OliveGreen,inline,backgroundcolor=OliveGreen!25,bordercolor=OliveGreen,#1]{#2}}
\newcommandx{\improvement}[2][1=]{\todo[linecolor=Plum,inline,backgroundcolor=Plum!25,bordercolor=Plum,#1]{#2}}

\usepackage[font=small,skip=5pt]{caption}

\newtheorem{thm}{Theorem}
\newtheorem{lem}{Lemma}

\newtheorem{cor}{Corollary}

\newtheorem{rem}{Remark}
\newtheorem{ass}{Assumption}

\newcolumntype{C}[1]{>{\centering\let\newline\\\arraybackslash\hspace{0pt}}m{#1}}

\newcommand\tagthis{\addtocounter{equation}{1}\tag{\theequation}}

\DeclareMathOperator{\R}{\mathbb{R}} 
\DeclareMathOperator{\Ocal}{\mathcal{O}} 
\newcommand{\eqdef}{\stackrel{\text{def}}{=}}
\newcommand{\setn}{[n]}
\renewcommand{\top}{T}

\makeatletter
\newcounter{subthm} 
\let\savedc@thm\c@hyp

\newcommand{\normhyp}{%
  \let\c@hyp\savedc@hyp 
  \renewcommand\thehyp{\arabic{hyp}}%
} 
\makeatother

\makeatletter
\newcounter{subass} 
\let\savedc@ass\c@hyp

\makeatother

%% file: sections/introduction.tex
\section{Introduction}\label{sec_intro}
We are interested in solving the \textit{finite-sum smooth} minimization problem
\begin{align*}
\min_{w \in \mathbb{R}^d} \left\{ F(w) = \frac{1}{n} \sum_{i=1}^n f_i(w) \right\}, \tagthis \label{emperical_risk_loss} 
\end{align*}
where each $f_i$, $i \in \setn\eqdef\{1,\dots,n\}$, has a Lipschitz continuous gradient. Throughout the paper, we consider the case where $F$ has a finite lower bound $F^*$. 

Problems of form \eqref{emperical_risk_loss} cover a wide range of convex and nonconvex problems in machine learning applications including but not limited to logistic regression, neural networks, multi-kernel learning, etc. In many of these applications, the number of component functions $n$ is very large, which makes the classical Gradient Descent (GD) method less efficient since it requires to compute a full gradient many times. 
Instead, a traditional alternative is to employ stochastic gradient descent (SGD)~\cite{RM1951,pegasos,bottou2016optimization}. 
In recent years, a large number  of improved variants of stochastic gradient algorithms called variance reduction methods have emerged, in particular, SAG/SAGA~\cite{SAGjournal,SAGA}, SDCA~\cite{SDCA}, MISO \cite{mairal2013optimization}, SVRG/S2GD~\cite{SVRG,S2GD}, SARAH~\cite{Nguyen2017_sarah}, etc. These methods were first analyzed for strongly convex problems of form \eqref{emperical_risk_loss}. Due to recent interest in deep neural networks, {\em nonconvex} problems of form \eqref{emperical_risk_loss} have been studied and analyzed by considering a number of different approaches including many variants of variance reduction techniques (see e.g. \cite{nonconvexSVRG,LihuaNonconvex,allen2017natasha,allen2017natasha2,fang2018spider}, etc.) 

We study the SARAH algorithm \cite{Nguyen2017_sarah,Nguyen2017_sarahnonconvex} depicted in  Algorithm~\ref{sarah}, slightly modified. We use upper index $s$ to indicate the $s$-th outer loop and lower index $t$ to indicate the $t$-th iteration in the inner loop.
The key update rule is
\begin{align*}
    v_{t}^{(s)} = \nabla f_{i_{t}} (w_{t}^{(s)}) - \nabla f_{i_{t}}(w_{t-1}^{(s)}) + v_{t-1}^{(s)}.
    \tagthis \label{sarah_update}
\end{align*}
The computed $v_{t}^{(s)}$ is used to update 
\begin{equation}
w_{t+1}^{(s)} = w_{t}^{(s)} - \eta v_{t}^{(s)}. \label{update_iter}
\end{equation}
 After $m$ iteration in the inner loop, the outer loop remembers the last computed $w_{m+1}^{(s)}$ and starts its loop anew -- first with a full gradient computation before again entering the inner loop with updates (\ref{sarah_update}).
Instead of remembering $\tilde{w}_s = w_{m+1}^{(s)}$ for the next outer loop, the original SARAH algorithm in \cite{Nguyen2017_sarah} uses $\tilde{w}_s = w_{t}^{(s)}$ with $t$ chosen uniformly at random from $\{0,1,\dots,m\}$. The authors of \cite{Nguyen2017_sarah} chose to do this in order to being able to analyze the convergence rate for a single outer loop -- since in practice it makes sense to keep the last computed $w_{m+1}^{(s)}$ if multiple outer loop iterations are used, we give full credit of Algorithm~\ref{sarah} to  \cite{Nguyen2017_sarah} and call this SARAH. We would like to attain an \textit{$\epsilon$-accurate solution} satisfying $\mathbb{E}[ \|\nabla F(\tilde{w})\|^2] \leq \epsilon$ for the outputted approximation $\tilde{w}$. 

\begin{algorithm}
   \caption{SARAH (modified of \cite{Nguyen2017_sarah})}
   \label{sarah}
\begin{algorithmic}
   \STATE {\bfseries Parameters:} the learning rate $\eta > 0$, the inner loop size $m$, and the outer loop size $S$
   \STATE {\bfseries Initialize:} $\tilde{w}_0$
   \STATE {\bfseries Iterate:}
   \FOR{$s=1,2,\dots,S $}
   \STATE $w_0^{(s)} = \tilde{w}_{s-1}$
   \STATE $v_0^{(s)} = \frac{1}{n}\sum_{i=1}^{n} \nabla f_i(w_0^{(s)})$
   \STATE $w_1^{(s)} = w_0^{(s)} - \eta v_0^{(s)}$
   \STATE {\bfseries Iterate:}
   \FOR{$t=1,\dots,m$}
   \STATE Sample $i_{t}$ uniformly at random from $\setn$
   \STATE $v_{t}^{(s)} = \nabla f_{i_{t}} (w_{t}^{(s)}) - \nabla f_{i_{t}}(w_{t-1}^{(s)}) + v_{t-1}^{(s)}$
   \STATE $w_{t+1}^{(s)} = w_{t}^{(s)} - \eta v_{t}^{(s)}$
   \ENDFOR
   \STATE Set $\tilde{w}_s = w_{m+1}^{(s)}$ \textbf{(modified point)}
   \ENDFOR
\end{algorithmic}
\end{algorithm} 

We will analyze SARAH for smooth nonconvex optimization, i.e., we study \eqref{emperical_risk_loss} with the following assumption
\begin{ass}[average-$L$-smooth]
\label{ass_avgLsmooth}
The objective function $F$ is $L$-average-smooth, i.e., there exists a constant $L > 0$ such that, $\forall w,w' \in \mathbb{R}^d$, 
\begin{gather*}
\mathbb{E}_{i} [\| \nabla f_i(w) - \nabla f_i(w') \|^2] = \frac{1}{n} \sum_{i=1}^n \| \nabla f_i(w) - \nabla f_i(w') \|^2 \leq L^2 \| w - w' \|^2. \tagthis\label{eq:Lsmooth_average}
\end{gather*} 
\end{ass}

We notice that, the above assumption is weaker than the assumption on $L$-smoothness of each $f_i$, $i=1,\dots,n$. Throughout this paper for non-convex results, we only consider Assumption \ref{ass_avgLsmooth} and no other assumptions. We stress that our convergence analysis only relies on the above average smooth assumption without bounded variance assumption (as required in \cite{LihuaNonconvex,zhou2018stochastic}).  

We note that Assumption \ref{ass_avgLsmooth} implies that $F$ is $L$-smooth, that is, there exists a constant $L > 0$ such that, $\forall w,w' \in \mathbb{R}^d$, $\| \nabla F(w) - \nabla F(w') \| \leq L \| w - w' \|$. By Theorem 2.1.5 in \cite{nesterov2004}, we could obtain
\begin{align*}
    F(w) \leq F(w') + \nabla F(w')^T(w-w') + \frac{L}{2}\|w-w'\|^2.
\tagthis\label{eq:Lsmooth} 
\end{align*}

We measure the convergence rate in terms of  total complexity $\mathcal{T}$, i.e., the total number of gradient computations. For SARAH we have
$$\mathcal{T}= S\cdot (n+2m).$$
We notice that SARAH, using the notation and definition of \cite{fang2018spider},  is a random algorithm ${\cal A}$ that maps functions $f$ to a sequence of iterates
$$[{\bf x}^k;i_k] = {\cal A}^{k}({\bf \xi}, \nabla f_{i_0}({\bf x}^0), \nabla f_{i_1}({\bf x}^1), \ldots, \nabla f_{i_{k-1}}({\bf x}^{k-1})),$$
where ${\cal A}^{k-1}$ is a measure mapping, $i_k$ is the individual function chosen by ${\cal A}$ at iteration $k$, and ${\bf \xi}$ is a uniform random vector with entries in $[0,1]$. Rephrasing Theorem~3 in  \cite{fang2018spider} states the following lower bound: For $n \leq \Ocal(\epsilon^{-2})$, there exists a function $f$ such that in order to find a point $\tilde{{\bf x}}$ for which accuracy $\|\nabla F(\tilde{{\bf x}})\|^2\leq \epsilon$, ${\cal A}$ must have a total complexity $\mathcal{T}$  of at least $\Omega(L\sqrt{n}/\epsilon)$ 
stochastic gradient computations. Applying this bound to SARAH tells us that if the final output $\tilde{w}_S$ has
$$ \mathbb{E}[\| \nabla F(\tilde{w}_S)\|^2]\leq \epsilon \mbox{ then } \mathcal{T}=S\cdot (n+2m) = \Omega(L\sqrt{n}/\epsilon), \ n \leq \Ocal(\epsilon^{-2}).$$
In this paper, we show that in SARAH we can choose parameters $S$ and $m$ such that, for $n \leq \Ocal(\epsilon^{-2})$, the total complexity is
$$ \mathcal{T}=S\cdot (n+2m) = {\cal O}(L\sqrt{n}/\epsilon).$$








\vspace{.25cm}

\noindent
{\bf Related Work:} The paper that introduces SARAH \cite{Nguyen2017_sarahnonconvex} is only able to analyze convergence of a single outer loop giving a total complexity of $\Ocal(n + \frac{L^2}{\epsilon^2})$. 

Besides the lower bound,  \cite{fang2018spider}  introduces SPIDER, as a variant of SARAH, which achieves the  best known convergence result in the nonconvex case.
SPIDER uses the SARAH update rule (\ref{sarah_update}) as was originally proposed in \cite{Nguyen2017_sarah} and the mini-batch version of SARAH in \cite{Nguyen2017_sarahnonconvex}. SPIDER and SARAH are different in terms of iteration (\ref{update_iter}), which are $w_{t+1} = w_{t} - \eta (v_t / \|v_t\|)$ and $w_{t+1} = w_{t} - \eta v_t$, respectively. Also, SPIDER does not divide into outer loop and inner loop as SARAH does although SPIDER does also perform a full gradient update after a certain fixed number of iterations.
A recent technical report \cite{wang2018spiderboost} provides an improved version of SPIDER called SpiderBoost which allows a larger learning rate. 
Both SPIDER and SpiderBoost are able to show for smooth nonconvex optimization a total complexity of $\Ocal\left( n + L \sqrt{n}/\epsilon \right)$. 


\begin{table*}
 \centering 
 \scriptsize
\caption{Comparison of results on the total complexity for smooth nonconvex optimization}
\label{table_nonconvex}
\begin{tabular}{|c|c|c|c| }
\hline 
\textbf{Method} & \textbf{Complexity} & \textbf{Additional assumption} \\
\hline
\hline
GD \cite{nesterov2004}  & $\mathcal{O}\left(\frac{n}{\epsilon} \right)$ & None \\
\hline
SVRG \cite{nonconvexSVRG}  & $\mathcal{O}\left( n + \frac{n^{2/3}}{\epsilon} \right)$ & None \\
\hline
SCSG  \cite{LihuaNonconvex}  & $\mathcal{O}\left( \left( \frac{\sigma}{\epsilon} \wedge n \right) + \frac{1}{\epsilon} \left( \frac{\sigma}{\epsilon} \wedge n \right)^{2/3} \right)$  & Bounded variance \\
& $\mathcal{O}\left( n + \frac{n^{2/3}}{\epsilon} \right)$  & None ($\sigma \to \infty$) \\
\hline
SNVRG \cite{zhou2018stochastic}  & $\mathcal{O}\left( \log^3\left( \frac{\sigma}{\epsilon} \wedge n \right) \left[ \left( \frac{\sigma}{\epsilon} \wedge n \right) + \frac{1}{\epsilon} \left( \frac{\sigma}{\epsilon} \wedge n \right)^{1/2} \right]  \right)$  & Bounded variance \\
& $\mathcal{O}\left( \log^3\left( n \right) \left( n + \frac{\sqrt{n}}{\epsilon} \right)  \right)$  & None ($\sigma \to \infty$) \\
\hline
SPIDER  \cite{fang2018spider}  & $\mathcal{O}\left( n + \frac{\sqrt{n}}{\epsilon} \right)$  & None \\
\hline
SpiderBoost  \cite{wang2018spiderboost}  & $\mathcal{O}\left( n + \frac{\sqrt{n}}{\epsilon} \right)$  & None \\
\hline
{\color{red}\textbf{SARAH (this paper)}} & $\textcolor{red}{\mathcal{O}\left(  \frac{\sqrt{n}}{\epsilon} \vee n  \right) }$ & {\color{red}\textbf{None}} \\
\hline
\end{tabular}
\end{table*}

Table~\ref{table_nonconvex}\footnote{$a \wedge b$ is defined as $\min\{a,b\}$ and $a \vee b$ is defined as $\max\{a,b\}$} shows the comparison of results on the total complexity for smooth nonconvex optimization. (a) Each of the complexities in Table~\ref{table_nonconvex} also depends on the Lipschitz constant $L$, however,  since we consider smooth optimization and it is custom to assume/design $L = \Ocal(1)$, we ignore the dependency on $L$ in the complexity results.  (b) Although many algorithms have appeared during the past few years, we only compare  algorithms having a convergence result which only supposes the smooth assumption.  
(c) Among algorithms with convergence results that only suppose the smooth assumption, Table~\ref{table_nonconvex} only mentions  recent state-of-the-art results. 
   (d) Although the bounded variance assumption $\mathbb{E}[\| \nabla f_i(w) - \nabla F(w) \|^2] \leq \sigma^2$ is acceptable in many existing literature, this additional assumption limits the applicability of these convergence results since it adds dependence on $\sigma$ which can be arbitrarily large. 
For fair comparison with  convergence analysis without the bounded variance assumption, $\sigma$ must be set to go to infinity -- and this is what is mentioned in Table~\ref{table_nonconvex}. 
As an example, from Table~\ref{table_nonconvex}  we observe that SCSG has an advantage over SVRG only if $\sigma = \Ocal(1)$ but, theoretically, it has the same total complexity as SVRG if $\sigma \to \infty$. 
(e) For completeness, incompatibility with assuming a bounded gradient $\mathbb{E}[\| \nabla f_i(w) \|^2] \leq \sigma$ has been discussed in \cite{Nguyen2018} for strongly convex objective functions.

According to the results in Table~\ref{table_nonconvex}, we can observe that SARAH enjoys the same fast convergence rate as those of SPIDER and SpiderBoost in the nonconvex case for finding a first-order stationary point based on only the average smooth assumption. Its complexity matches the lower-bound worst case complexity in \cite{fang2018spider} up to a constant factor when $n \leq \Ocal(\epsilon^{-2})$. 

\vspace{.25cm}

\noindent
{\bf Contributions:}
We summarize our key contributions as follows. 

\begin{itemize}
    \item {\em Smooth Non-Convex.} 
    We provide a convergence analysis for the full SARAH algorithm with 
multiple outer iterations for nonconvex problems (unlike in \cite{Nguyen2017_sarahnonconvex} which  only analyses a single outer iteration). Its complexity matches the lower-bound worst case complexity in \cite{fang2018spider} up to a constant factor when $n \leq \Ocal(\epsilon^{-2})$. 
The convergence analysis \textbf{only} supposes the average smooth assumption (which is weaker than Lipschitz continuous assumption on each component gradient) in the non-convex case 
 (Theorem~\ref{thm_main_nonconvex}). 
 We extend this result to the \textit{mini-batch} case (Theorem~\ref{thm_main_nonconvex_mb}). 
 \item {\em Smooth Convex.} In order to complete the picture, we study SARAH+ \cite{Nguyen2017_sarah} which was designed as a variant of SARAH for convex optimization.
We propose a novel variant of SARAH+ called \textit{SARAH++}. 
Here, we study the {\em iteration complexity} measured by the total number of iterations (which counts one full gradient computation as adding one iteration to the complexity) -- and leave an analysis of the total complexity as an open problem. 
For SARAH++ we show
a sublinear convergence rate in the general convex case (Theorem~\ref{thm_general_convex_sarahpp}) and a linear convergence rate in the strongly convex case (Theorem~\ref{thm_strongly_convex_sarahpp}). SARAH itself may already lead to good convergence and there may no need to introduce SARAH++; in numerical experiments we show the advantage of SARAH++ over SARAH. We further propose a practical version called \textit{SARAH Adaptive} which improves the performance of SARAH and SARAH++ for convex problems -- numerical experiments on various data sets show good overall performance. 
\item For the convergence analysis of SARAH for the non-convex case and SARAH++ for the convex case we show that the analysis generalizes the total complexity of Gradient Descent (GD) (Remarks~\ref{rem_main_derivation_nonconvex_2} and \ref{rem_sarahpp}), i.e., the analysis reproduces known total complexity results of GD. Up to the best of our knowledge, this is the first variance reduction method having this property. 
\end{itemize}

%% file: sections/method.tex
\section{Non-Convex Case: Convergence Analysis of SARAH}\label{sec_analysis}

SARAH is very different from other algorithms since it has a \textit{biased} estimator of the gradient. Therefore, in order to analyze SARAH's convergence rate,  it is non-trivial to use existing proof techniques  from unbiased estimator algorithms such as SGD, SAGA, and SVRG.





\subsection{A single batch case}

We start analyzing SARAH (Algorithm \ref{sarah}) for the case where we choose a single sample $i_t$ uniformly at random from $[n]$ in the inner loop. 

\begin{lem}\label{lem_main_derivation_nonconvex}
Suppose that Assumption \ref{ass_avgLsmooth} holds. Consider a single outer loop iteration in SARAH (Algorithm~\ref{sarah})  with $\eta \leq \frac{2}{L(\sqrt{1 + 4m} + 1)}$. Then, for any $s \geq 1$, we have
\begin{equation}
\mathbb{E}[ F(w^{(s)}_{m+1}) ]  \leq \mathbb{E} [ F(w^{(s)}_0) ] - \frac{\eta}{2} \sum_{t=0}^{m} \mathbb{E}[ \| \nabla F(w^{(s)}_{t})\|^2 ].   \label{eq:001}
\end{equation} 
\end{lem}

The above result is for a single outer loop iteration of SARAH, which includes a full gradient step together with the inner loop. 
Since the outer loop iteration concludes with $\tilde{w}_s = w^{(s)}_{m+1}$, and $\tilde{w}_{s-1} = w^{(s)}_{0}$,
we have
\begin{align*}
    \mathbb{E}[ F(\tilde{w}_s) ]  &\leq \mathbb{E} [ F(\tilde{w}_{s-1}) ] - \frac{\eta}{2} \sum_{t=0}^{m} \mathbb{E}[ \| \nabla F(w_{t}^{(s)})\|^2 ]. 
\end{align*}
Summing over $1\leq s\leq S$ gives
\begin{align*}
    \mathbb{E}[ F(\tilde{w}_S) ] & \leq \mathbb{E}[ F(\tilde{w}_0) ] - \frac{\eta}{2} \sum_{s=1}^S \sum_{t=0}^{m} \mathbb{E}[ \| \nabla F(w_{t}^{(s)})\|^2 ]. \tagthis \label{eq_recursive}
\end{align*}
This proves our main result:

\begin{thm}[Smooth nonconvex]\label{thm_main_nonconvex}
Suppose that Assumption \ref{ass_avgLsmooth} holds. Consider SARAH (Algorithm~\ref{sarah}) with $\eta \leq \frac{2}{L(\sqrt{1 + 4m} + 1)}$. Then, for any given $\tilde{w}_0$, we have
\begin{align*}
    & \frac{1}{(m+1)S}\sum_{s=1}^S \sum_{t=0}^{m} \mathbb{E}[ \| \nabla F(w_{t}^{(s)})\|^2 ]  \leq \frac{2}{\eta [(m+1) S]} [ F(\tilde{w}_0) - F^*],
\end{align*}
where $F^*$ is any lower bound of $F$, and $w_{t}^{(s)}$ is the result of the $t$-th iteration in the $s$-th outer loop. 
\end{thm} 
The proof easily follows from \eqref{eq_recursive} since $F^*$ is a lower bound of $F$ (that is, $\mathbb{E}[ F(\tilde{w}_S) ] \geq F^*$). We note that the term
\begin{align*}
    \frac{1}{(m+1)S}\sum_{s=1}^S \sum_{t=0}^{m} \mathbb{E}[ \| \nabla F(w_{t}^{(s)})\|^2 ]
\end{align*}
is simply the average of the expectation of the squared norms of the gradients of all the iteration results generated by SARAH. For nonconvex problems, our goal is to achieve
\begin{align*}
    \frac{1}{(m+1)S}\sum_{s=1}^S \sum_{t=0}^{m} \mathbb{E}[ \| \nabla F(w_{t}^{(s)})\|^2 ] \leq \epsilon. 
\end{align*}


We note that, for simplicity, if $\bar{w}_s$ is chosen uniformly at random from all the iterations generated by SARAH, we are able to have accuracy $\mathbb{E}[ \| \nabla F(\bar{w}_s)\|^2 ] \leq \epsilon$. 

\begin{cor}\label{cor_main_derivation_nonconvex}
Suppose that Assumption \ref{ass_avgLsmooth} holds. Consider SARAH (Algorithm~\ref{sarah}) with $\eta = \Ocal(\frac{1}{L \sqrt{m+1}})$ where $m$ is the inner loop size. Then, in order to achieve an $\epsilon$-accurate solution, the total complexity is $$\Ocal\left(  \left[ \left( \frac{n + 2m}{\sqrt{m+1}} \right) \frac{1}{\epsilon} \right] \vee \left[ n + 2m \right] \right).$$ 
\end{cor}

The total complexity can be minimized over the inner loop size $m$. 
By choosing $m = n$, we achieve the minimal total complexity:
\begin{cor}\label{cor_main_derivation_nonconvex_2}
Suppose that Assumption \ref{ass_avgLsmooth} holds. Consider SARAH (Algorithm~\ref{sarah}) with $\eta = \Ocal(\frac{1}{L \sqrt{m+1}})$ where $m$ is the inner loop size and chosen equal to $m = n$. Then, in order to achieve an $\epsilon$-accurate solution, the total complexity is 
$$\Ocal\left( \frac{\sqrt{n}}{\epsilon} \vee n \right).$$ 
\end{cor}

\begin{rem}\label{rem_main_derivation_nonconvex_2}
The total complexity in Corollary~\ref{cor_main_derivation_nonconvex} covers all 
choices for the inner loop size $m$. For example, in the case of $m=0$, SARAH recovers the Gradient Descent (GD) algorithm which has total complexity $\Ocal\left( \frac{n}{\epsilon} \right)$. Theorem~\ref{thm_main_nonconvex} for $m=0$ also recovers the requirement on the learning rate for GD, which is $\eta \leq \frac{1}{L}$. 
\end{rem}

The above results explain the relationship between SARAH and GD  and explains the advantages of the inner loop and outer loop of SARAH. SARAH becomes more beneficial in ML applications where $n$ is large.

\subsection{Mini-batch case}

The above results can be extended to the \textit{mini-batch} case where instead of choosing a single sample $i_t$, we choose $b$ samples uniformly at random from $[n]$ for updating $v_t$ in the inner loop. We then replace $v_t$ in Algorithm \ref{sarah} by
\begin{align*}
    v_{t}^{(s)} = \frac{1}{b} \sum_{i \in I_{t}} [\nabla f_{i} (w_{t}^{(s)}) - \nabla f_{i}(w_{t-1}^{(s)})] + v_{t-1}^{(s)}, \tagthis \label{sarah_update_mb}
\end{align*}
where we choose a mini-batch $I_{t} \subseteq \setn$ of size $b$ uniformly at random at each iteration of the inner loop. The result of Theorem~\ref{thm_main_nonconvex} generalizes as follows. 
\begin{thm}[Smooth nonconvex with mini-batch]\label{thm_main_nonconvex_mb}
Suppose that Assumption \ref{ass_avgLsmooth} holds. Consider SARAH (Algorithm~\ref{sarah}) by replacing $v_t$ in the inner loop size by \eqref{sarah_update_mb} with 
\begin{align*}
    \eta \leq \frac{2}{L\left(\sqrt{1 + \frac{4m}{b}\left(\frac{n-b}{n-1} \right)} + 1\right)}. 
\end{align*}
Then, for any given $\tilde{w}_0$, we have
\begin{align*}
    & \frac{1}{(m+1)S}\sum_{s=1}^S \sum_{t=0}^{m} \mathbb{E}[ \| \nabla F(w_{t}^{(s)})\|^2 ] \leq \frac{2}{\eta [(m+1) S]} [ F(\tilde{w}_0) - F^*],
\end{align*}
where $F^*$ is any lower bound of $F$, and $w_{t}^{(s)}$ is the $t$-th iteration in the $s$-th outer loop.  
\end{thm} 

We can again derive similar corollaries as was done for Theorem~\ref{thm_main_nonconvex}.  

\begin{cor}\label{cor_main_derivation_nonconvex_mb}
For the conditions in Theorem~\ref{thm_main_nonconvex_mb}, in order to achieve an $\epsilon$-accurate solution, the total complexity is $$\Ocal\left(  \left[ \left( \frac{n + 2 b m}{m + 1}    \right)  \left( \sqrt{1 + \frac{4m}{b}\left(\frac{n-b}{n-1} \right)} \right) \frac{1}{\epsilon} \right] \vee \left[ n + 2 b m \right] \right).$$ 
\end{cor}

\begin{cor}\label{cor_main_derivation_nonconvex_2_mb}
For the conditions in Theorem~\ref{thm_main_nonconvex_mb} and Corollary~\ref{cor_main_derivation_nonconvex_mb} with $b = n^{\alpha}$ and $m = n^{\beta}$ where $\alpha + \beta = 1$ with $\beta \geq 1/2$ and $0 \leq \alpha \leq 1/2$, in order to achieve an $\epsilon$-accurate solution, the total complexity is $$\Ocal\left( \frac{\sqrt{n}}{\epsilon} \vee n \right).$$ 
\end{cor}


\section{Convex Case: SARAH++: A New Variant of SARAH+}\label{sec_sarahpp}

In this section, we propose a new variant of SARAH+ (Algorithm~\ref{sarah_plus}) \cite{Nguyen2017_sarah}, called SARAH++ (Algorithm~\ref{sarah_new}), for convex problems of form \eqref{emperical_risk_loss}. 

Different from SARAH, SARAH+ provides a stopping criteria for the inner loop; as soon as 
$$\|v_{t-1}^{(s)}\|^2 \leq  \gamma \|v_{0}^{(s)}\|^2,$$ the inner loop finishes. This idea originates from the property of SARAH that, for each outer loop iteration $s$,  $\mathbb{E}[\|v_t^{(s)}\|^2] \to 0$ as $t \to \infty$ in the strongly convex case (Theorems 1a and 1b in \cite{Nguyen2017_sarah}). Therefore, it does not make any sense to update with tiny steps when $\|v_{t}^{(s)}\|^2$ is small. (We note that SVRG \cite{SVRG} does not have this property.) SARAH+ suggests to empirically choose parameter $\gamma = 1/8$ \cite{Nguyen2017_sarah} without theoretical guarantee. 
\begin{algorithm} 
   \caption{SARAH+ \cite{Nguyen2017_sarah}}
   \label{sarah_plus}
\begin{algorithmic}
   \STATE {\bfseries Parameters:} the learning rate $\eta > 0$, $0 < \gamma \leq 1$, the maximum inner loop size $m$, and the outer loop size $S$ 
   \STATE {\bfseries Initialize:} $\tilde{w}_0$
   \STATE {\bfseries Iterate:}
   \FOR{$s=1,2,\dots,S$}
   \STATE $w_0^{(s)} = \tilde{w}_{s-1}$
   \STATE $v_0^{(s)} = \frac{1}{n}\sum_{i=1}^{n} \nabla f_i(w_0^{(s)})$
   \STATE $w_1^{(s)} = w_0^{(s)} - \eta v_0^{(s)}$
   \STATE $t = 1$
   \WHILE{$\|v_{t-1}^{(s)}\|^2 >  \gamma \|v_{0}^{(s)}\|^2$ {\bf and} $t \leq m$}
   \STATE Sample $i_{t}$ uniformly at random from $\setn$
   \STATE $v_{t}^{(s)} = \nabla f_{i_{t}} (w_{t}^{(s)}) - \nabla f_{i_{t}}(w_{t-1}^{(s)}) + v_{t-1}^{(s)}$
   \STATE $w_{t+1}^{(s)} = w_{t}^{(s)} - \eta v_{t}^{(s)}$
   \STATE $t \leftarrow t + 1$
   \ENDWHILE
   \STATE Set $\tilde{w}_s = w_{t}^{(s)}$
   \ENDFOR
\end{algorithmic}
\end{algorithm} 

Here, we modify SARAH+ (Algorithm~\ref{sarah_plus}) into SARAH++ (Algorithm~\ref{sarah_new}) by choosing the stopping criteria for the inner loop as 
$$\|v_{t-1}^{(s)}\|^2 <  \gamma \|v_{0}^{(s)}\|^2 \mbox{ where } \gamma \geq L\eta$$ 
and by introducing a stopping criteria for the outer loop.


\subsection{Details SARAH++ and Convergence Analysis}

Before analyzing and explaining SARAH++ in detail, we introduce the following assumptions used in this section. 

\begin{ass}[$L$-smooth]
\label{ass_Lsmooth}
Each $f_i: \mathbb{R}^d \to \mathbb{R}$, $i \in \setn$, is $L$-smooth, i.e., there exists a constant $L > 0$ such that, $\forall w,w' \in \mathbb{R}^d$, 
\begin{gather*}
\| \nabla f_i(w) - \nabla f_i(w') \| \leq L \| w - w' \|. \tagthis\label{eq:Lsmooth_basic}
\end{gather*} 
\end{ass}

\begin{ass}[$\mu$-strongly convex]
\label{ass_stronglyconvex}
The function $F: \mathbb{R}^d \to \mathbb{R}$, is $\mu$-strongly convex, i.e., there exists a constant $\mu > 0$ such that $\forall w,w' \in \mathbb{R}^d$, 
\begin{gather*}
F(w)   \geq  F(w') + \nabla F(w')^\top (w - w') + \tfrac{\mu}{2}\|w - w'\|^2.
\end{gather*}
\end{ass}

Under Assumption \ref{ass_stronglyconvex}, let us define the (unique) optimal solution of \eqref{emperical_risk_loss} as $w_{*}$.
Then strong convexity of $F$ implies that 
 \begin{equation} \label{eq:strongconvexity2}
 2\mu [ F(w) - F(w_{*})] \leq  \| \nabla F(w)\|^2, \ \forall w \in \mathbb{R}^d. 
 \end{equation}
We note here, for future use, that for strongly convex functions of the form \eqref{emperical_risk_loss}, arising in machine learning applications, the condition number is defined as $\kappa\eqdef L/\mu$. 
Assumption \ref{ass_stronglyconvex} covers a wide range of problems, e.g. $l_2$-regularized empirical risk minimization problems with convex losses. 

We separately assume the special case of strong convexity of all $f_i$'s with $\mu=0$, called the general convexity assumption, which we will use for convergence analysis.
\begin{ass}
\label{ass_convex}
Each function $f_i: \mathbb{R}^d \to \mathbb{R}$, $i \in \setn$, is convex, i.e.,
\begin{gather*}
f_i(w)    \geq f_i(w') + \nabla f_i(w')^\top (w - w'). 
\end{gather*}
\end{ass}

SARAH++ is motivated by the following lemma. 
\begin{lem}\label{lem_basic_lem_01}
Suppose that Assumptions \ref{ass_Lsmooth} and \ref{ass_convex} hold. Consider a single outer loop iteration in SARAH (Algorithm~\ref{sarah})  with $\eta \leq \frac{1}{L}$. Then, for $t \geq 0$ and any $s \geq 1$, we have
\begin{align*}
   & \mathbb{E}[F(w^{(s)}_{t+1}) - F(w_*)] \leq \mathbb{E}[F(w^{(s)}_{t}) - F(w_*)]  - 
   \frac{\eta}{2} \mathbb{E} [ \| \nabla F(w^{(s)}_t) \|^2]  + 
   \frac{\eta}{2} \left( L\eta \mathbb{E}[\| v^{(s)}_0 \|^2] - \mathbb{E}[\| v^{(s)}_t \|^2] \right), \tagthis \label{eq_lem_01}
\end{align*}
where $w_*$ is any optimal solution of $F$.
\end{lem}



Clearly, if 
$$L\eta \mathbb{E}[\| v_0^{(s)} \|^2] - \mathbb{E}[\| v_t^{(s)} \|^2] \leq \gamma \mathbb{E}[\| v_0^{(s)} \|^2] - \mathbb{E}[\| v_t^{(s)} \|^2] \leq 0,$$
where $\eta \leq \frac{\gamma}{L}$, inequality \eqref{eq_lem_01} implies
\begin{align*}
    \mathbb{E}[F(w^{(s)}_{t+1}) - F(w_*)] & \leq \mathbb{E}[F(w^{(s)}_{t}) - F(w_*)]   - \frac{\eta}{2} \mathbb{E} [ \| \nabla F(w^{(s)}_t) \|^2]. 
\end{align*}
For this reason, we choose the stopping criteria for the inner loop in SARAH++ as 
$\| v_t^{(s)} \|^2 < \gamma \|v_0^{(s)} \|^2$ with $\gamma\geq L\eta$. 
Unlike SARAH+, for analyzing the convergence rate $\gamma$ can be as small as $L\eta$.

\begin{algorithm}[t]
   \caption{SARAH++}
   \label{sarah_new}
\begin{algorithmic}
   \STATE {\bfseries Parameters:} The controlled factor $0< \gamma \leq 1$, the learning rate $0 < \eta \leq \frac{\gamma}{L}$, the total iteration $T > 0$, and the maximum inner loop size $m \leq T$. 
   \STATE {\bfseries Initialize:} $\tilde{w}_0$
   \STATE $G = 0$
   \STATE {\bfseries Iterate:}
   \STATE $s = 0$
   \WHILE{$G < T$}
   \STATE $s \leftarrow s + 1$
   \STATE $w_0^{(s)} = \tilde{w}_{s-1}$
   \STATE $v_0^{(s)} = \frac{1}{n}\sum_{i=1}^{n} \nabla f_i(w_0^{(s)})$
   \STATE $t = 0$
   \WHILE{$\|v_{t}^{(s)}\|^2 \geq \gamma \|v_{0}^{(s)}\|^2$ {\bf and} $t \leq m$}
   \STATE $w_{t+1}^{(s)} = w_{t}^{(s)} - \eta v_{t}^{(s)}$
   \STATE $t \leftarrow t + 1$
   \IF{$m \neq 0$}
   \STATE Sample $i_{t}$ uniformly at random from $\setn$
   \STATE $v_{t}^{(s)} = \nabla f_{i_{t}} (w_{t}^{(s)}) - \nabla f_{i_{t}}(w_{t-1}^{(s)}) + v_{t-1}^{(s)}$
   \ENDIF
   \ENDWHILE
   \STATE $T_s = t$
   \STATE $\tilde{w}_s = w_{T_s}^{(s)}$
   \STATE $G \leftarrow G + T_s$
   \ENDWHILE
   \STATE $S = s$
   \STATE Set $\hat{w} = \tilde{w}_S$
\end{algorithmic}
\end{algorithm}

The above discussion leads to SARAH++ (Algorithm~\ref{sarah_new}). In order to analyze its convergence for convex problems, we define random variable $T_s$ as the stopping time of the inner loop in the $s$-th outer iteration: 
\begin{align*}
    T_s = \min \left\{ \min_{t \geq 0} \left\{ t : \| v_{t}^{(s)} \|^2 < \gamma \| v_0^{(s)} \|^2 \right\} , m + 1 \right\} \ , \ s = 1,2,\dots
\end{align*}
Note that $T_s$ is at least 1 since at $t = 0$, the condition $\| v_{0}^{(s)} \|^2 \geq \gamma \| v_0^{(s)} \|^2$ always holds (and $m \geq 0$). 

Let random variable $S$  be the stopping time of the outer iterations as a function of an algorithm parameter $T > 0$: 
\begin{align*}
    S = \min_{\hat{S}} \left\{ \hat{S} : \sum_{s=1}^{\hat{S}} T_s \geq T  \right\}. 
\end{align*}
Notice that SARAH++ maintains a running sum $G=\sum_{j=1}^s T_i$ against which parameter $T$ is compared in the stopping criteria of the outer loop.

For the general convex case which supposes Assumption~\ref{ass_convex} in addition to smoothness we have the next theorem.

\begin{thm}[Smooth general convex]\label{thm_general_convex_sarahpp}
Suppose that Assumptions \ref{ass_Lsmooth} and \ref{ass_convex} hold. Consider SARAH++ (Algorithm~\ref{sarah_new}) with $\eta \leq \frac{\gamma}{L}$, $0 < \gamma \leq 1$. Then, 
$$\mathbb{E}\left[ \frac{1}{T_1 + \dots + T_S}  \sum_{s=1}^S \sum_{t=0}^{T_s-1} \mathbb{E} [ \| \nabla F(w_t^{(s)}) \|^2 |T_1,\dots,T_S ]  \right],$$
the expectation of the average of the squared norm of the gradients of all iterations generated by SARAH++, is bounded by 
\begin{align*}
    \frac{2}{T \eta} [ F(\tilde{w}_0) - F(w_*) ].
\end{align*}
\end{thm}

The theorem leads to the next corollary about iteration complexity, i.e.,  we bound $T$ which is the total number of iterations performed by the inner loop across all outer loop iterations. 
This is different from the total complexity since $T$ does not separately count the $n$ gradient evaluations when the full gradient is computed in the outer loop.

\begin{cor}[Smooth general convex]\label{cor_general_convex_sarahpp}
For the conditions in Theorem~\ref{thm_general_convex_sarahpp} with $\eta = \Ocal(\frac{1}{L})$, we achieve an $\epsilon$-accurate solution after $\Ocal(\frac{1}{\epsilon})$ inner loop iterations.
\end{cor}

By supposing Assumption~\ref{ass_stronglyconvex} in addition to the smoothness and general convexity assumptions, we can prove a linear convergence rate. For strongly convex objective functions we have the following  result.
\begin{thm}[Smooth strongly convex]\label{thm_strongly_convex_sarahpp}
Suppose that Assumptions \ref{ass_Lsmooth}, \ref{ass_stronglyconvex} and \ref{ass_convex} hold. Consider SARAH++ (Algorithm~\ref{sarah_new}) with $\eta \leq \frac{\gamma}{L}$, $0 < \gamma \leq 1$. Then, for the final output $\hat{w}$ of SARAH++, we have
\begin{align*}
    \mathbb{E} [ F(\hat{w}) - F(w_*) ] &\leq (1 - \mu \eta)^T [ F(\tilde{w}_0) - F(w_*) ]. \tagthis \label{eq_thm_strongly_convex}
\end{align*}
\end{thm}



This leads to the following iteration complexity. 
\begin{cor}[Smooth strongly convex]\label{cor_strongly_convex_sarahpp}
For the conditions in Theorem~\ref{thm_strongly_convex_sarahpp} with $\eta = \Ocal(\frac{1}{L})$, we achieve $\mathbb{E} [ F(\hat{w}) - F(w_*) ] \leq \epsilon$ after $\Ocal(\kappa \log(\frac{1}{\epsilon}))$ total iterations, where $\kappa = L/\mu$ is the condition number.
\end{cor}

\begin{rem}\label{rem_sarahpp}
The proofs of the above results hold for any $m\leq T$. 
If we choose $m=0$, then SARAH++ reduces to the Gradient Descent algorithm since the inner ``while'' loop stops right after updating $w_1^{(s)} = w_0^{(s)} - \eta v_0^{(s)}$. 
In this case, Corollaries~\ref{cor_general_convex_sarahpp} and \ref{cor_strongly_convex_sarahpp} recover the rate of convergence and complexity of GD. 
\end{rem}

In  this section, we showed that \textit{SARAH++ has a guarantee of theoretical convergence} (see Theorems~\ref{thm_general_convex_sarahpp} and \ref{thm_strongly_convex_sarahpp}) while SARAH+ does not have such a guarantee.



An interesting open question we would like to discuss here is the total complexity of SARAH++. Although we have shown the convergence results of SARAH++ in terms of the iteration complexity, the total complexity which is computed as the total number of evaluations of the component gradient functions still remains an open question. It is clear  that the total complexity must depend on the learning rate $\eta$ (or $\gamma$) -- the factor that decides when to stop the inner iterations. 

We note that $T$ can be ``closely'' understood as the total number of  updates $w^{(s)}_{t+1}$ of the algorithm. The total complexity is equal to 
$\sum_{i=1}^S (n + 2 (T_i - 1))$. 
For the special case $T_i = 1$, $i = 1,\dots,S$, the algorithm recovers the GD algorithm with $T = \sum_{i=1}^S T_s = S$. Since each full gradient takes $n$ gradient evaluations, the total complexity for this case is equal to $n S = \Ocal(\frac{n}{\epsilon})$ (in the general convex case) and $n S = \Ocal(n \kappa \log(\frac{1}{\epsilon}))$ (in the strongly convex case). 


However, it is non-trivial to derive the total complexity of SARAH++ since it should depend on the learning rate $\eta$. We leave this question as an open direction for future research.

%% file: sections/sarahpp_experiments.tex
\subsection{Numerical Experiments}\label{sec_experiments}

 

Paper \cite{Nguyen2017_sarah} provides experiments showing good overall performance of SARAH over other algorithms such as SGD \cite{RM1951}, SAG \cite{SAG}, SVRG \cite{SVRG}, etc. For this reason, we provide experiments comparing SARAH++ directly with SARAH. We notice that SARAH (with multiple outer loops) like SARAH++ has theoretical guarantees with sublinear convergence for general convex and linear convergence for strongly convex problems as proved in \cite{Nguyen2017_sarah}. Because of these theoretical guarantees (which SARAH+ does not have), SARAH itself may already perform well for convex problems and the question is whether SARAH++ offers an advantage.
 

We consider $\ell_2$-regularized logistic regression problems with
\begin{align*}
f_i(w) = \log(1 + \exp(-y_i \langle x_i, w\rangle )) +
\tfrac{\lambda}{2} \| w \|^2, \tagthis \label{eq_logistic_regression}
\end{align*}
where $\{x_i,y_i\}_{i=1}^n$ is the training data and the regularization parameter $\lambda$ is set to $1/n$, a widely-used
value in  literature \cite{SAG,Nguyen2017_sarah}. The condition number is equal to $\kappa = L/\mu = n$. We conducted experiments to demonstrate the advantage in performance of SARAH++ over SARAH for convex problems on popular data sets including \textit{covtype} ($n = 406,708$ training data; estimated $L \simeq 1.90$) and \textit{ijcnn1} ($n = 91,701$ training data; estimated $L \simeq 1.77$) from LIBSVM \cite{LIBSVM}. 
\begin{figure}[h]
 \centering
 \includegraphics[width=0.36\textwidth]{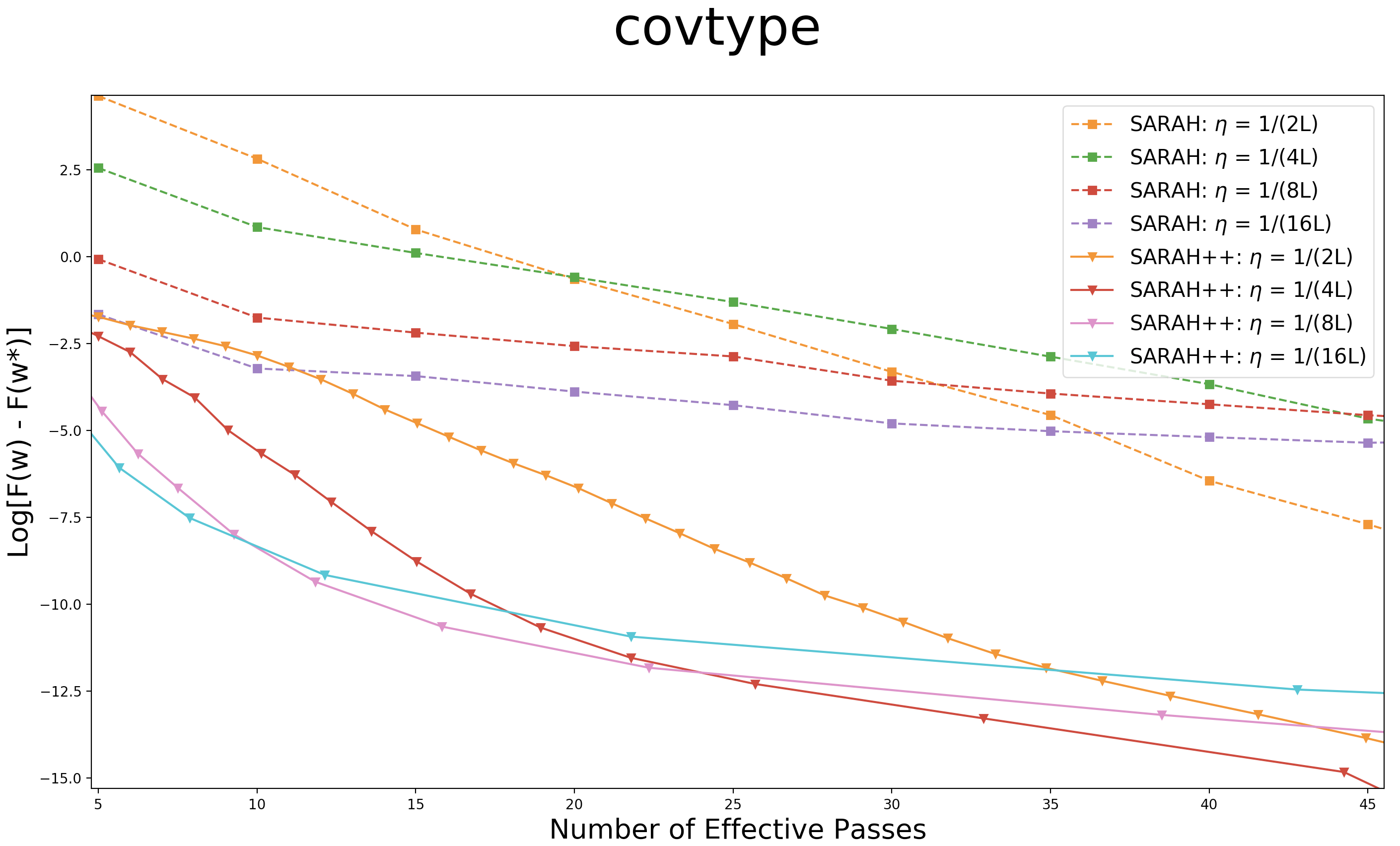}
 \includegraphics[width=0.36\textwidth]{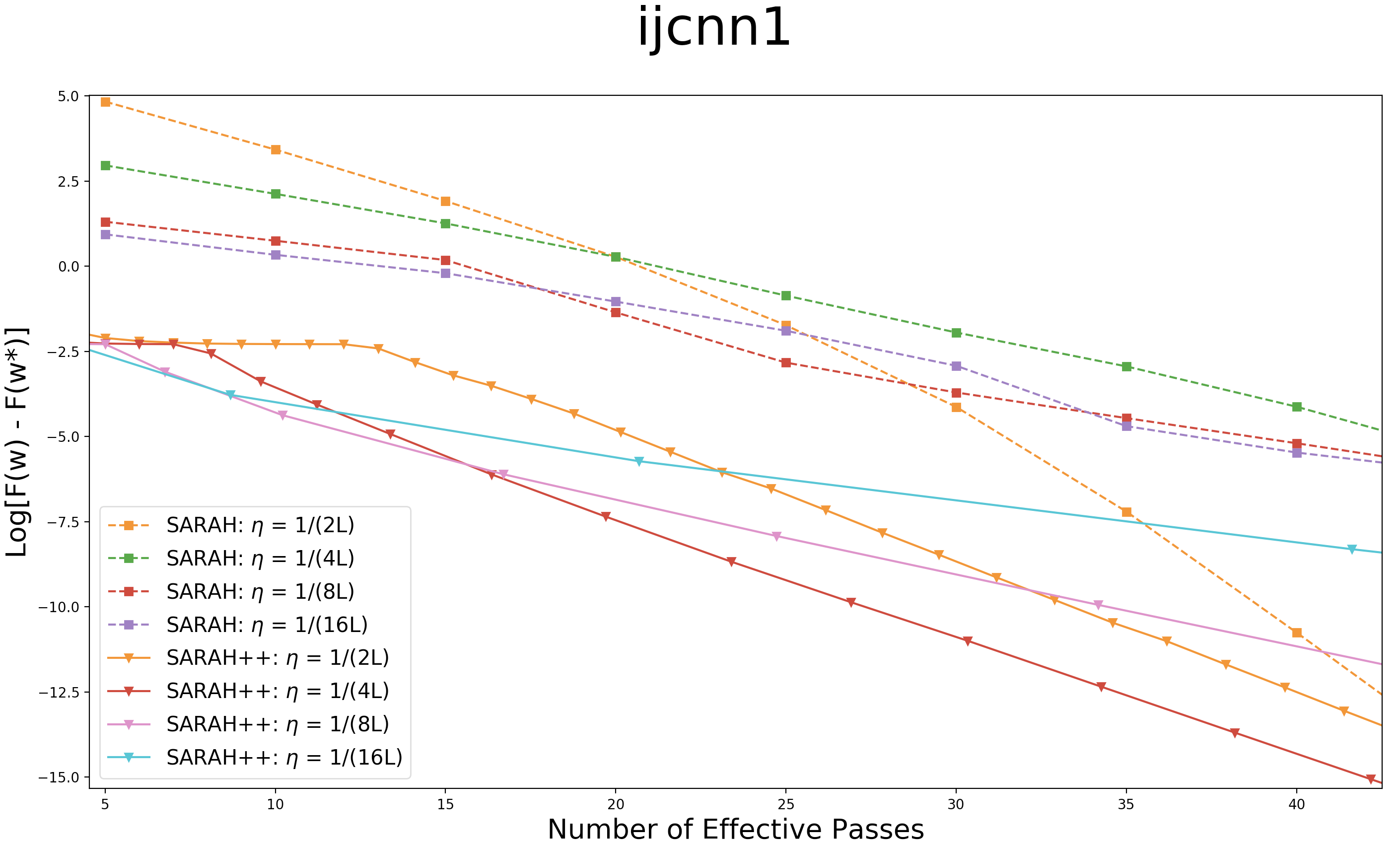}
  \caption{Comparisons of $\log[F(w) - F(w_*)]$ between SARAH++ and SARAH with different learning rates on \textit{covtype} and \textit{ijcnn1} datasets}
  \label{fig_sarahpp_sarah}
 \end{figure} 
 
Figure~\ref{fig_sarahpp_sarah} shows comparisons between SARAH++ and SARAH for different values of learning rate $\eta$. We depicted the value of $\log [F(w) - F(w_*)] $ (i.e. $F(w) - F(w_*)$ in log scale) for the $y$-axis and ``number of effective passes'' (or number of epochs, where an epoch is the equivalent of $n$ component gradient evaluations or one full gradient computation) for the $x$-axis. For SARAH, we choose the outer loop size $S = 10$ and tune the inner loop size $m = \{0.5n, n, 2n, 3n, 4n \}$ to achieve the best performance. The optimal solution $w_*$ of the strongly convex problem in \eqref{eq_logistic_regression} is found by using Gradient Descent with stopping criterion $\| \nabla F(w) \|^2 \leq 10^{-15}$. We observe that, SARAH++ achieves improved overall performance compared to regular SARAH as shown in Figure~\ref{fig_sarahpp_sarah}. From the experiments we see that the stopping criteria $\| v_t^{(s)} \|^2 < \gamma \|v_0^{(s)}\|^2$ ($\gamma = L\eta$) of SARAH++ is indeed important. 
The stopping criteria helps the inner loop to prevent updating tiny redundant steps. 

%% file: sections/adaptive.tex
\subsection{SARAH Adaptive: A New Practical Variant}\label{sec_adaptive}

We now propose a practical adaptive method 
which aims to improve performance. Although we do not have any theoretical result for this adaptive method, numerical experiments are very promising and they heuristically show the improved performance on different data sets.

\begin{algorithm}
   \caption{SARAH Adaptive}
   \label{sarah_adap}
\begin{algorithmic}
   \STATE {\bfseries Parameters:} The maximum inner loop size $m$, and the outer loop size $S$, the factor $0 < \gamma \leq 1$. 
   \STATE {\bfseries Initialize:} $\tilde{w}_0$
   \STATE {\bfseries Iterate:}
   \FOR{$s=1,2,\dots,S $}
   \STATE $w_0^{(s)} = \tilde{w}_{s-1}$
   \STATE $v_0^{(s)} = \frac{1}{n}\sum_{i=1}^{n} \nabla f_i(w_0^{(s)})$
   \STATE $t=0$
   \WHILE{$\|v_{t}^{(s)}\|^2 \geq \gamma \|v_{0}^{(s)}\|^2$ {\bf and} $t \leq m$}
   \STATE $\eta_t^{(s)} = \frac{1}{L} \cdot \frac{\| v_{t}^{(s)} \|^2}{\| v_0^{(s)} \|^2}$ \textbf{(adaptive)}
   \STATE $w_{t+1}^{(s)} = w_{t}^{(s)} - \eta_t^{(s)} v_{t}^{(s)}$
   \STATE $t \leftarrow t+1$
   \IF{$m \neq 0$}
   \STATE Sample $i_{t}$ uniformly at random from $\setn$
   \STATE $v_{t}^{(s)} = \nabla f_{i_{t}} (w_{t}^{(s)}) - \nabla f_{i_{t}}(w_{t-1}^{(s)}) + v_{t-1}^{(s)}$
   \ENDIF
   \ENDWHILE
   \STATE Set $\tilde{w}_s = w_{t}^{(s)}$
   \ENDFOR
\end{algorithmic}
\end{algorithm} 

The motivation of this algorithm comes from the intuition of Lemma~\ref{lem_basic_lem_01} (for convex optimization). For a single outer loop with $\eta \leq \frac{1}{L}$, \eqref{eq_lem_01} holds for SARAH (Algorithm~\ref{sarah}). Hence, for any $s$, we intentionally choose $\eta = \eta_t^{(s)} = \frac{\| v_{t}^{(s)} \|^2}{L \| v_0^{(s)} \|^2}$ such that $L\eta \mathbb{E}[\| v_0^{(s)} \|^2] - \mathbb{E}[\| v_t^{(s)} \|^2] = 0$. Since $\| v_{t}^{(s)} \|^2 \leq \| v_0^{(s)} \|^2$, $t \geq 0$, in \cite{Nguyen2017_sarah} for convex problems, we have $\eta_t^{(s)} \leq \frac{1}{L}$, $t \geq 0$. We also stop the inner loop by the stopping criteria $\|v_{t}^{(s)}\|^2 < \gamma \|v_{0}^{(s)}\|^2$ for some $0 < \gamma \leq 1$. SARAH Adaptive is given in detail in Algorithm~\ref{sarah_adap} without convergence analysis. 
 

\begin{figure}[h]
 \centering
 \includegraphics[width=0.36\textwidth]{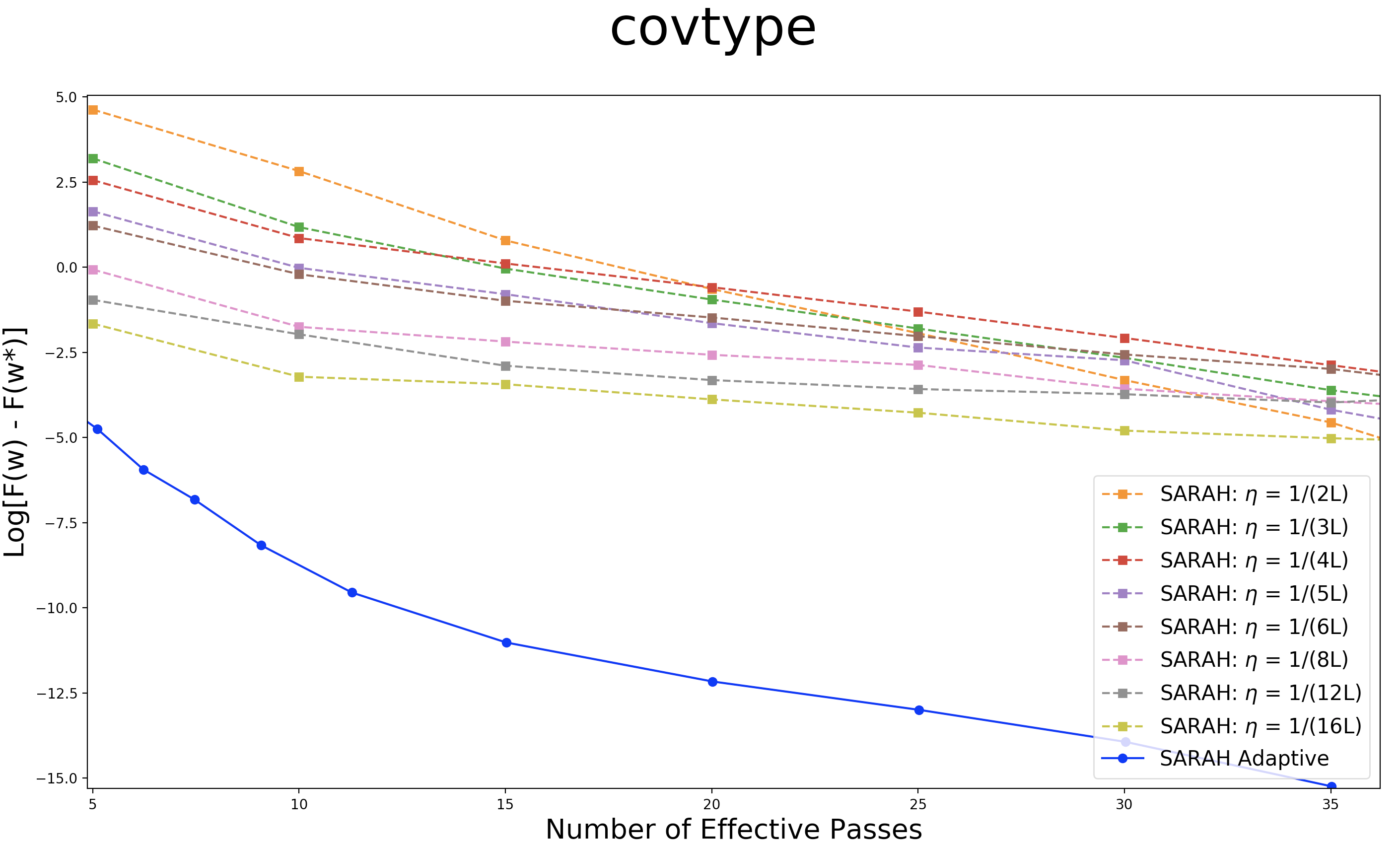}
 \includegraphics[width=0.36\textwidth]{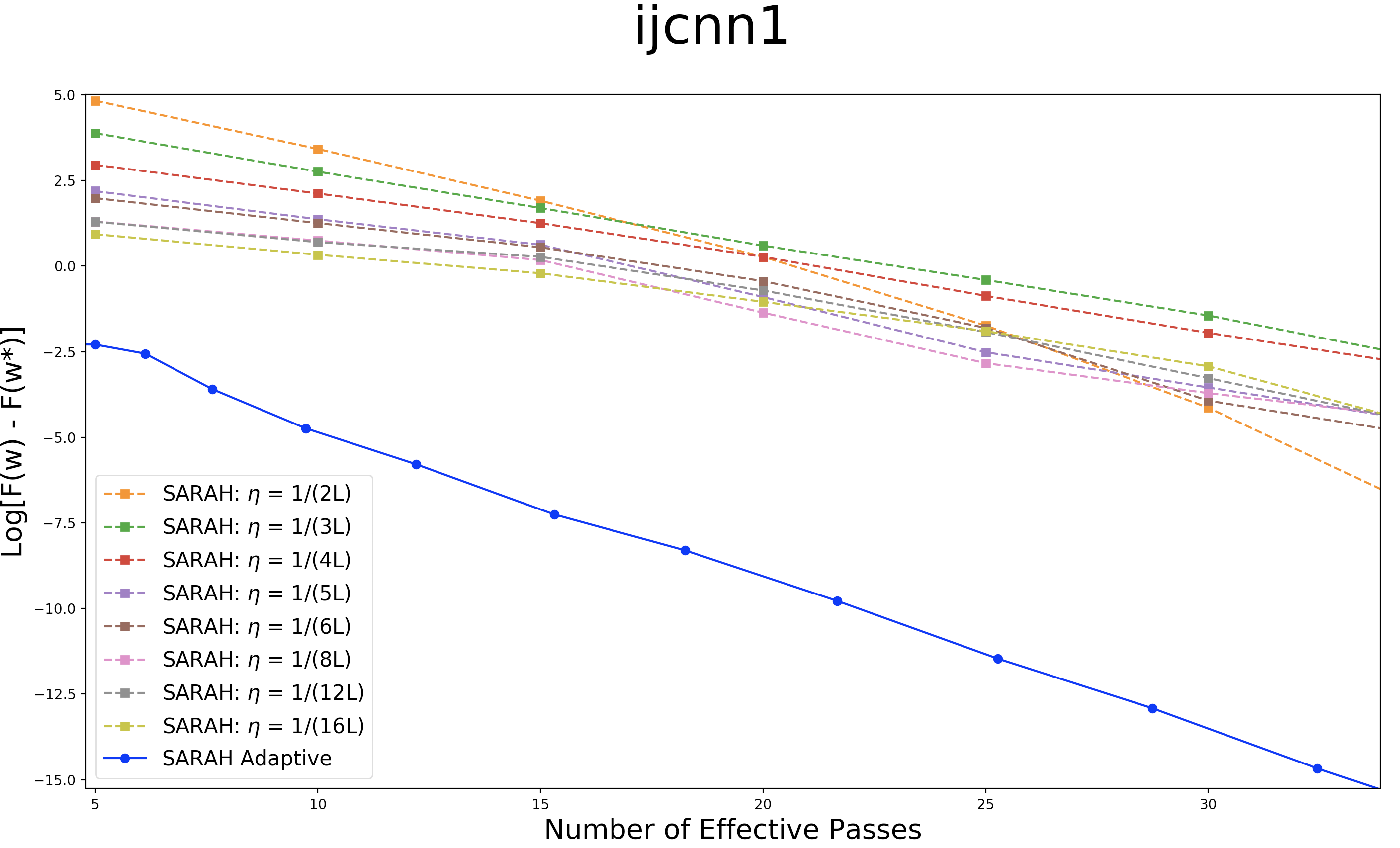}
  \caption{Comparisons of $\log[F(w) - F(w_*)]$ between SARAH Adaptive and SARAH with different learning rates on \textit{covtype} and \textit{ijcnn1} datasets}
  \label{fig_sarah_adap_sarah}
 \end{figure}

\begin{figure}[h!]
 \centering
 \includegraphics[width=0.36\textwidth]{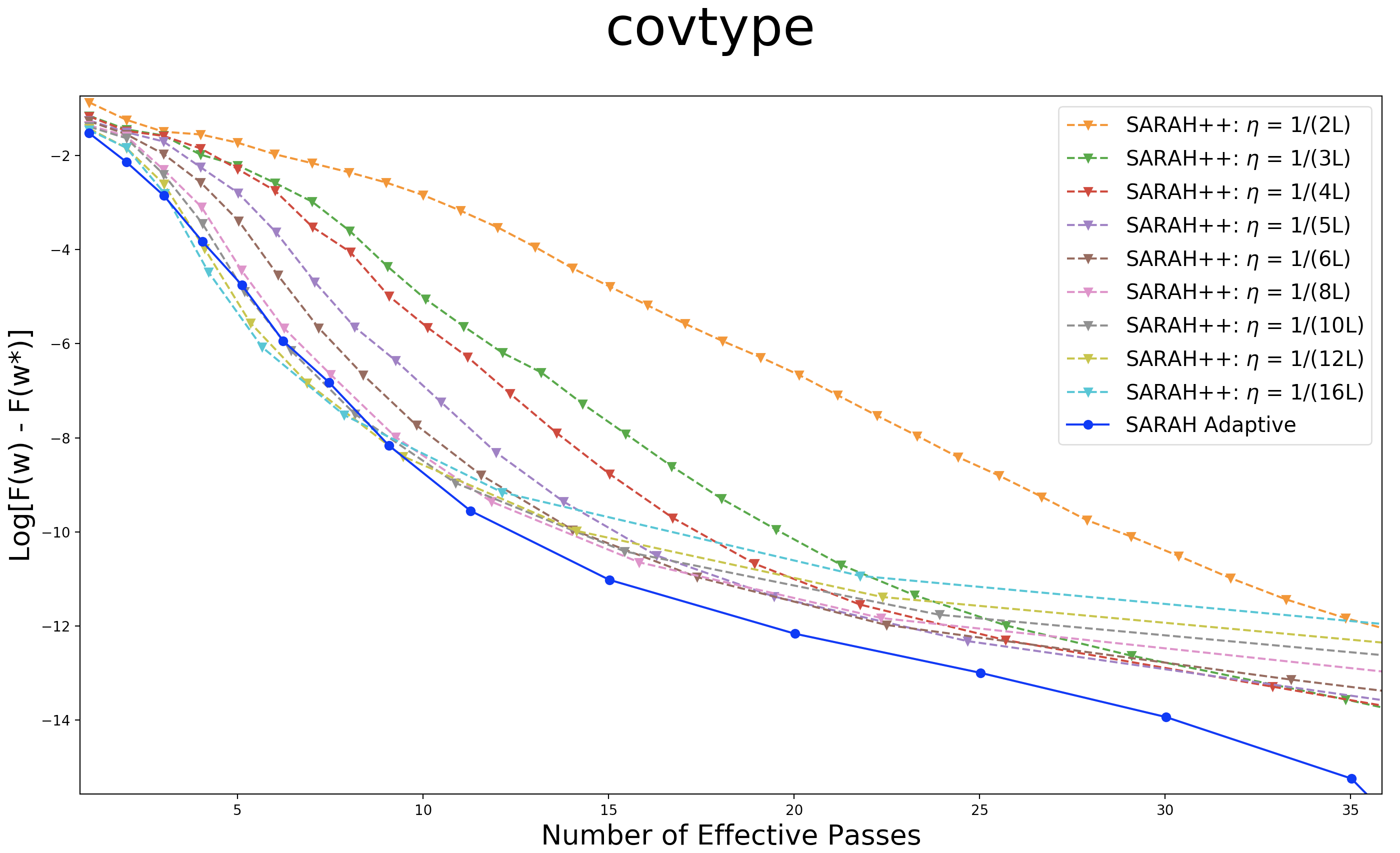}
 \includegraphics[width=0.36\textwidth]{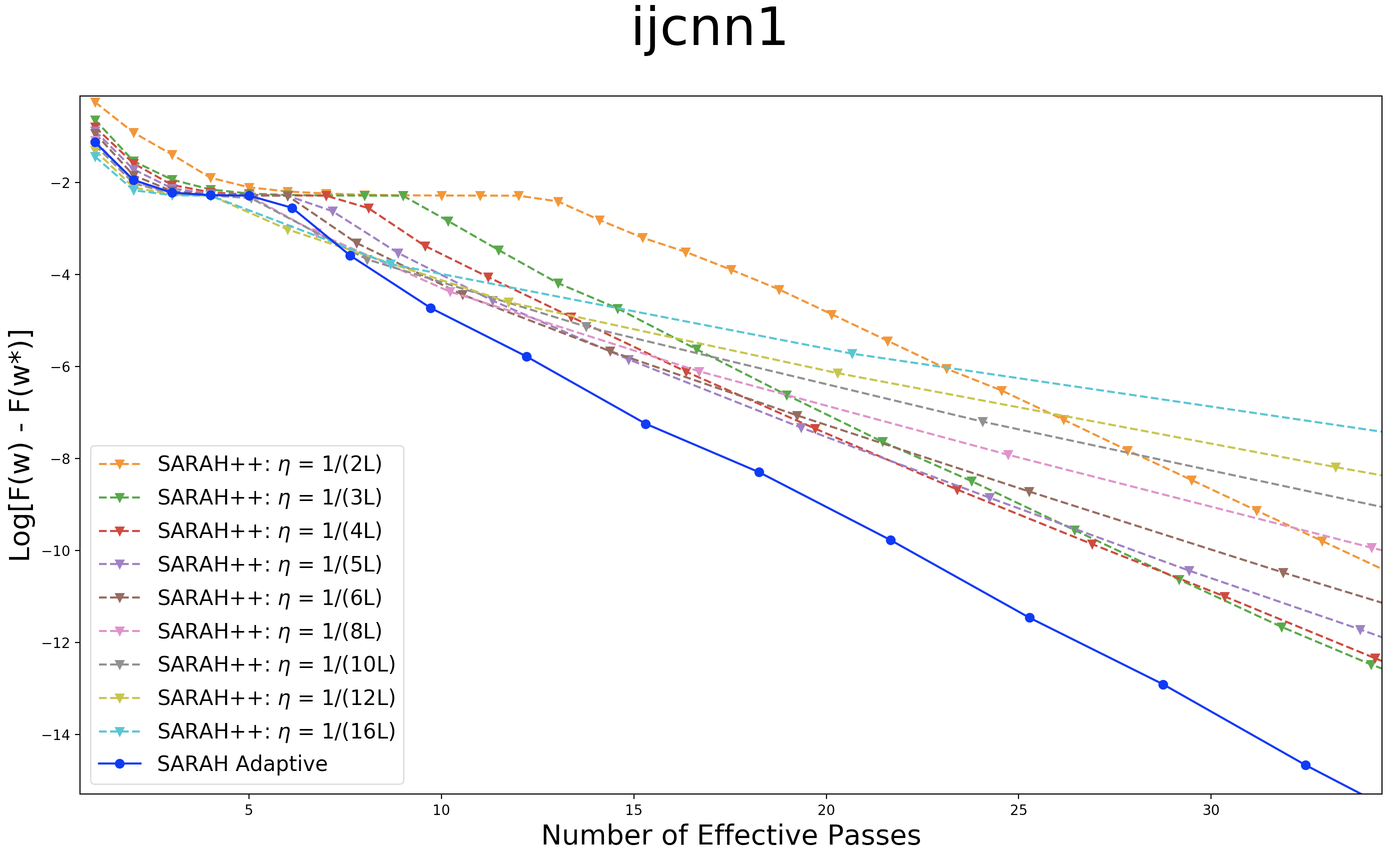}
  \caption{Comparisons of $\log[F(w) - F(w_*)]$ between SARAH Adaptive and SARAH++ with different learning rates on \textit{covtype} and \textit{ijcnn1} datasets}
  \label{fig_sarah_adap_sarahpp}
 \end{figure}
 
We have conducted numerical experiments on the same datasets and problems as introduced in the previous subsection. 
Figures~\ref{fig_sarah_adap_sarah} and \ref{fig_sarah_adap_sarahpp} show the comparison between SARAH Adaptive and SARAH and SARAH++ for different values of $\eta$. We observe that SARAH Adaptive has an improved performance over SARAH and SARAH++ (without tuning learning rate). We also present the numerical performance of SARAH Adaptive for different values of $\gamma$ in Appendix. We also present the numerical performance of SARAH Adaptive for different values of $\gamma$ in Appendix. 

We note that additional experiments in this section on more data sets are performed in Appendix. 

%% file: sections/conclusion.tex




%% file: sections/appendix.tex
\section*{Appendix}

\section*{Useful Existing Results}


\begin{lem}[Lemma 2 in \cite{Nguyen2017_sarah} (or in \cite{Nguyen2017_sarahnonconvex})]\label{lem:var_diff_01}
Consider $v_{t}^{(s)}$ defined by \eqref{sarah_update} (or \eqref{sarah_update_mb}) in SARAH (Algorithm \ref{sarah}) for any $s \geq 1$. Then for any $t\geq 1$, 
\begin{align*}
&\mathbb{E}[ \| \nabla F(w_{t}^{(s)}) - v_{t}^{(s)} \|^2 ] 
= \sum_{j = 1}^{t} \mathbb{E}[ \| v_{j}^{(s)} - v_{j-1}^{(s)} \|^2 ]  - \sum_{j = 1}^{t} \mathbb{E}[ \| \nabla F(w_{j}^{(s)}) - \nabla F(w_{j-1}^{(s)}) \|^2 ]. \tagthis\label{eq_lemma_sarah}
\end{align*}
\end{lem}

\begin{lem}[Lemma 3 in \cite{Nguyen2017_sarah}]\label{lem_bound_var_diff_str_02}
Suppose that Assumptions \ref{ass_Lsmooth} and \ref{ass_convex} hold. Consider $v_{t}^{(s)}$ defined as \eqref{sarah_update} in SARAH (Algorithm \ref{sarah}) with $\eta < 2/L$ for any $s \geq 1$. Then we have that for any $t\geq 0$, 
\begin{align*}
\mathbb{E}[ \| \nabla F(w_{t}^{(s)}) - v_{t}^{(s)} \|^2 ] 
\leq  \frac{\eta L}{2 - \eta L} \Big[ \mathbb{E}[ \|v_{0}^{(s)} \|^2] - \mathbb{E}[\| v_{t}^{(s)} \|^2] \Big].
\tagthis\label{eq:bound1}
\end{align*}
\end{lem}

\section*{Nonconvex SARAH}

\subsection*{Proof of Lemma \ref{lem_main_derivation_nonconvex}}

\textbf{Lemma \ref{lem_main_derivation_nonconvex}}.
\textit{Suppose that Assumption \ref{ass_avgLsmooth} holds. Consider SARAH (Algorithm~\ref{sarah}) within a single outer loop with $\eta \leq \frac{2}{L(\sqrt{1 + 4m} + 1)}$. Then, for any $s \geq 1$, we have}
\begin{align*}
\mathbb{E}[ F(w_{m+1}^{(s)}) ]  &\leq \mathbb{E} [ F(w_0^{(s)}) ] - \frac{\eta}{2} \sum_{t=0}^{m} \mathbb{E}[ \| \nabla F(w_{t}^{(s)})\|^2 ]. 
\end{align*} 

\begin{proof}
We use some parts of the proof in \cite{Nguyen2017_sarahnonconvex}. By Assumption \ref{ass_avgLsmooth} and $w_{t+1}^{(s)} = w_{t}^{(s)} - \eta v_{t}^{(s)}$, for any $s \geq 1$, we have
\begin{align*}
\mathbb{E}[ F(w_{t+1}^{(s)})] & \overset{\eqref{eq:Lsmooth}}{\leq}  \mathbb{E}[ F(w_{t}^{(s)})] - \eta \mathbb{E}[\nabla F(w_{t}^{(s)})^\top v_{t}^{(s)}] 
+ \frac{L\eta^2}{2} \mathbb{E} [ \| v_{t}^{(s)} \|^2 ] 
\\
& = \mathbb{E}[ F(w_{t}^{(s)})] - \frac{\eta}{2} \mathbb{E}[ \| \nabla F(w_{t}^{(s)})\|^2 ] 
+ \frac{\eta}{2} \mathbb{E}[ \| \nabla F(w_{t}^{(s)}) - v_{t}^{(s)} \|^2 ] 
- \left( \frac{\eta}{2} - \frac{L\eta^2}{2} \right) \mathbb{E} [ \| v_{t}^{(s)} \|^2 ], \tagthis \label{eq_proof_lem_000}
\end{align*}
where the last equality follows from the fact
$a^Tb = \frac{1}{2}\left[\|a\|^2 + \|b\|^2 - \|a-b\|^2\right],$ for any $a,b\in \R^d$. By summing over $t = 0,\dots,m$, we have
\begin{align*}
  \mathbb{E}[ F(w_{m+1}^{(s)}) ]  &\leq \mathbb{E} [ F(w_0^{(s)}) ] - \frac{\eta}{2} \sum_{t=0}^{m} \mathbb{E}[ \| \nabla F(w_{t}^{(s)})\|^2 ] \\ & \qquad + \frac{\eta}{2} \left( \sum_{t=0}^{m} \mathbb{E}[ \| \nabla F(w_{t}^{(s)}) - v_{t}^{(s)} \|^2 ]  
 - ( 1 - L\eta ) \sum_{t=0}^{m} \mathbb{E} [ \| v_{t}^{(s)} \|^2 ] \right). \tagthis \label{eq_proof_lem_001}
\end{align*}
Now, we would like to determine $\eta$ such that the expression in \eqref{eq_proof_lem_001}
\begin{align*}
    \sum_{t=0}^{m} \mathbb{E}[ \| \nabla F(w_{t}^{(s)}) - v_{t}^{(s)} \|^2 ]  - ( 1 - L\eta ) \sum_{t=0}^{m} \mathbb{E} [ \| v_{t}^{(s)} \|^2 ] \leq 0. 
\end{align*}
Let $\mathcal{F}_j = \sigma(w_0,w_1,\dots,w_j)$ be the $\sigma$-algebra generated by $w_0,w_1,\dots,w_j$. Note that $\mathcal{F}_j$ also contains all information of $v_0,\dots,v_{j-1}$. We have
\begin{align*}
\mathbb{E}[ \| v_{j}^{(s)} - v_{j-1}^{(s)} \|^2 | \mathcal{F}_j ] &\overset{\eqref{sarah_update}}{=} \mathbb{E}[ \| \nabla f_{i_{j}} (w_{j}^{(s)}) - \nabla f_{i_{j}}(w_{j-1}^{(s)}) \|^2 | \mathcal{F}_j ] \\ &
\overset{\eqref{eq:Lsmooth_average}}{\leq} L^2 \| w_{j}^{(s)} - w_{j-1}^{(s)} \|^2 = L^2 \eta^2 \| v_{j-1}^{(s)} \|^2, \ j \geq 1.
\end{align*}
Taking the expectations to both sides yields
\begin{align*}
    \mathbb{E}[ \| v_{j}^{(s)} - v_{j-1}^{(s)} \|^2] \leq L^2 \eta^2 \mathbb{E}[ \| v_{j-1}^{(s)} \|^2]. \tagthis \label{eq:afsag242}
\end{align*}
Hence, by Lemma \ref{lem:var_diff_01}, we have
\begin{align*}
\mathbb{E}[ \| \nabla F(w_{t}^{(s)}) - v_{t}^{(s)} \|^2 ] & \leq \sum_{j = 1}^{t} \mathbb{E}[ \| v_{j}^{(s)} - v_{j-1}^{(s)} \|^2 ] \overset{\eqref{eq:afsag242}}{\leq} L^2 \eta^2 \sum_{j = 1}^{t} \mathbb{E}[ \| v_{j - 1}^{(s)} \|^2 ]. 
\end{align*}
Note that $\| \nabla F(w_{0}^{(s)}) - v_{0}^{(s)} \|^2 = 0$. Hence, by summing over $t = 0,\dots,m$ ($m \geq 1$), we have
\begin{align*}
\sum_{t=0}^{m} \mathbb{E}\| \nabla F(w_{t}^{(s)}) - v_{t}^{(s)} \|^2 \leq L^2 \eta^2 \Big[ m  \mathbb{E}\|v_{0}^{(s)} \|^2 + (m-1) \mathbb{E}\|v_{1}^{(s)} \|^2 + \dots + \mathbb{E}\|v_{m-1}^{(s)}\|^2 \Big ]. 
\end{align*}
By choosing $\eta \leq \frac{2}{L\left(\sqrt{1 + 4m} + 1\right)}$, we have
\begin{align*}
& \sum_{t=0}^{m} \mathbb{E}[ \| \nabla F(w_{t}^{(s)}) - v_{t}^{(s)} \|^2 ]  - ( 1 - L\eta ) \sum_{t=0}^{m} \mathbb{E} [ \| v_{t}^{(s)} \|^2 ] \\
& \leq L^2 \eta^2 \Big[ m  \mathbb{E}\|v_{0}^{(s)} \|^2 + (m-1) \mathbb{E}\|v_{1}^{(s)} \|^2 + \dots + \mathbb{E}\|v_{m-1}^{(s)}\|^2 \Big ] - (1 - L\eta) \Big [ \mathbb{E}\|v_{0}^{(s)} \|^2 + \mathbb{E}\|v_{1}^{(s)} \|^2 + \dots + \mathbb{E}\|v_{m}^{(s)}\|^2  \Big ] \\
& \leq \Big[L^2\eta^2 m - (1 - L\eta) \Big] \sum_{t=1}^{m} \mathbb{E} [ \| v_{t-1}^{(s)} \|^2 ] \leq 0,  \tagthis \label{eq:equal_zero}
\end{align*}
since $\eta = \frac{2}{L\left(\sqrt{1 + 4m} + 1\right)}$ is a root of equation $L^2\eta^2 m - (1 - L\eta) = 0$. Therefore, with $\eta \leq \frac{2}{L(\sqrt{1 + 4m} + 1)}$, we have
\begin{align*}
\mathbb{E}[ F(w_{m+1}^{(s)}) ]  &\leq \mathbb{E} [ F(w_0^{(s)}) ] - \frac{\eta}{2} \sum_{t=0}^{m} \mathbb{E}[ \| \nabla F(w_{t}^{(s)})\|^2 ]. 
\end{align*}

\end{proof}


\subsection*{Proof of Corollary \ref{cor_main_derivation_nonconvex}}

\textbf{Corollary \ref{cor_main_derivation_nonconvex}. }
\textit{Suppose that Assumption \ref{ass_avgLsmooth} holds. Consider SARAH (Algorithm~\ref{sarah}) with $\eta = \Ocal(\frac{1}{L \sqrt{m+1}})$ where $m$ is the inner loop size. Then, in order to achieve $\epsilon$-accurate solution, the total complexity is $\Ocal\left( \left[ \left( \frac{n + 2m}{\sqrt{m+1}} \right) \frac{1}{\epsilon} \right] \vee \left[ n + 2m \right] \right)$.} 

\begin{proof}
In order to achieve 
\begin{align*}
    \frac{1}{(m+1)S}\sum_{s=1}^S \sum_{t=0}^{m} \mathbb{E}[ \| \nabla F(w_{t}^{(s)})\|^2 ] \leq \epsilon, 
\end{align*} 
we need
\begin{align*}
    \frac{2}{\eta [(m+1) S]} [ F(\tilde{w}_0) - F^*] = \epsilon. \tagthis \label{eq_cor_001}
\end{align*}
Let us choose $\eta$ such that
\begin{align*}
    \eta = \frac{2}{L(3\sqrt{m+1})} \overset{m\geq 0}{\leq} \frac{2}{L(2\sqrt{m+1}+1)} = \frac{2}{L(\sqrt{4m + 4} + 1)} \leq \frac{2}{L(\sqrt{1 + 4m} + 1)}. \tagthis \label{eq_cor_002}
\end{align*}
Hence, in order to achieve \eqref{eq_cor_001}, we need
\begin{align*}
    S &=  \frac{2}{\eta [(m+1) \epsilon]} [ F(\tilde{w}_0) - F^*] \overset{\eqref{eq_cor_002}}{=} \frac{3 L [ F(\tilde{w}_0) - F^*]}{(\sqrt{m+1})} \frac{1}{\epsilon} = \Ocal\left( \left[ \frac{1}{\sqrt{m+1}} \cdot \frac{1}{\epsilon} \right] \vee 1 \right) 
\end{align*}
since $S \geq 1$. Therefore, the total complexity to achieve $\epsilon$-accurate solution is $$(n + 2m)S = \Ocal\left( \left[ \left( \frac{n + 2m}{\sqrt{m+1}} \right) \frac{1}{\epsilon} \right] \vee \left[ n + 2m \right] \right).$$ 
\end{proof}


\subsection*{Proof of Theorem \ref{thm_main_nonconvex_mb}}

\textbf{Theorem \ref{thm_main_nonconvex_mb}} (Smooth nonconvex with mini-batch). 
\textit{Suppose that Assumption \ref{ass_avgLsmooth} holds. Consider SARAH (Algorithm~\ref{sarah}) by replacing $v_t$ in the inner loop size by \eqref{sarah_update_mb} with 
\begin{align*}
    \eta \leq \frac{2}{L\left(\sqrt{1 + \frac{4m}{b}\left(\frac{n-b}{n-1} \right)} + 1\right)}. 
\end{align*}
Then, for any given $\tilde{w}_0$, we have
\begin{align*}
    \frac{1}{(m+1)S}\sum_{s=1}^S \sum_{t=0}^{m} \mathbb{E}[ \| \nabla F(w_{t}^{(s)})\|^2 ]  \leq \frac{2}{\eta [(m+1) S]} [ F(\tilde{w}_0) - F^*],
\end{align*}
where $F^*$ is any lower bound of $F$, and $w_{t}^{(s)}$ is the $t$-th iteration in the $s$-th outer loop.  }

\begin{proof}
Following the proof of Lemma~\ref{lem_main_derivation_nonconvex}, we would like to determine $\eta$ such that the expression in \eqref{eq_proof_lem_001}
\begin{align*}
    \sum_{t=0}^{m} \mathbb{E}[ \| \nabla F(w_{t}^{(s)}) - v_{t}^{(s)} \|^2 ]  - ( 1 - L\eta ) \sum_{t=0}^{m} \mathbb{E} [ \| v_{t}^{(s)} \|^2 ] \leq 0. 
\end{align*}
Let
\begin{align*}
\xi_t = \nabla f_{t} (w_{j}^{(s)}) - \nabla f_{t}(w_{j-1}^{(s)}). \tagthis \label{xi_value}
\end{align*}

Let $\mathcal{F}_{j} = \sigma(w_0^{(s)},I_1,I_2,\dots,I_{j-1})$ be the $\sigma$-algebra generated by $w_0^{(s)},I_1,I_2,\dots,I_{j-1}$; $\mathcal{F}_{0} = \mathcal{F}_{1} = \sigma(w_0^{(s)})$. Note that $\mathcal{F}_{j}$ also contains all the information of $w_0^{(s)},\dots,w_{j}^{(s)}$ as well as $v_0^{(s)},\dots,v_{j-1}^{(s)}$. We have
\allowdisplaybreaks
\begin{align*}
& \mathbb{E}[\| v_j^{(s)} - v_{j-1}^{(s)} \|^2 | \mathcal{F}_{j} ] - \| \nabla F(w_{j}^{(s)}) - \nabla F(w_{j-1}^{(s)}) \|^2 \\
& \qquad \overset{\eqref{sarah_update_mb}}{=} \mathbb{E} \Big[ \Big\| \frac{1}{b} \sum_{i \in I_{j}} [\nabla f_{i} (w_{j}^{(s)}) - \nabla f_{i}(w_{j-1}^{(s)})] \Big\|^2 \Big| \mathcal{F}_{j} \Big] - \Big \| \frac{1}{n} \sum_{i=1}^{n} [\nabla f_{i} (w_{j}^{(s)}) - \nabla f_{i}(w_{j-1}^{(s)})] \Big \|^2 \\
& \qquad \overset{\eqref{xi_value}}{=} \mathbb{E} \Big[ \Big\| \frac{1}{b} \sum_{i \in I_{j}} \xi_i \Big\|^2 \Big| \mathcal{F}_{j} \Big] - \Big \| \frac{1}{n} \sum_{i=1}^{n} \xi_i \Big \|^2 \\
& \qquad = \frac{1}{b^2} \mathbb{E} \Big[ \sum_{i \in I_{j}} \sum_{k \in I_{j}} \xi_i^\top \xi_k  \Big| \mathcal{F}_{j} \Big] - \frac{1}{n^2} \sum_{i=1}^{n} \sum_{k=1}^{n} \xi_i^\top \xi_k \\
& \qquad = \frac{1}{b^2} \mathbb{E} \Big[ \sum_{i \neq k \in I_{j}} \xi_i^\top \xi_k + \sum_{i \in I_{j}} \xi_i^\top \xi_i  \Big| \mathcal{F}_{j} \Big] - \frac{1}{n^2} \sum_{i=1}^{n} \sum_{k=1}^{n} \xi_i^\top \xi_k \\
& \qquad = \frac{1}{b^2} \Big[ \frac{b}{n} \frac{(b-1)}{(n-1)} \sum_{i \neq k } \xi_i^\top \xi_k + \frac{b}{n}\sum_{i=1}^{n} \xi_i^\top \xi_i  \Big] - \frac{1}{n^2} \sum_{i=1}^{n} \sum_{k=1}^{n} \xi_i^\top \xi_k \\
& \qquad = \frac{1}{b^2} \Big[ \frac{b}{n} \frac{(b-1)}{(n-1)} \sum_{i=1}^{n} \sum_{k=1}^{n} \xi_i^\top \xi_k + \left( \frac{b}{n} - \frac{b}{n} \frac{(b-1)}{(n-1)} \right) \sum_{i=1}^{n} \xi_i^\top \xi_i  \Big] - \frac{1}{n^2} \sum_{i=1}^{n} \sum_{k=1}^{n} \xi_i^\top \xi_k \\
& \qquad = \frac{1}{b n} \Big[ \left( \frac{(b-1)}{(n-1)} - \frac{b}{n} \right) \sum_{i=1}^{n} \sum_{k=1}^{n} \xi_i^\top \xi_k + \frac{(n-b)}{(n-1)}  \sum_{i=1}^{n} \xi_i^\top \xi_i  \Big] \\
& \qquad = \frac{1}{b n} \left( \frac{n-b}{n-1} \right) \Big[ - \frac{1}{n}  \sum_{i=1}^{n} \sum_{k=1}^{n} \xi_i^\top \xi_k +  \sum_{i=1}^{n} \xi_i^\top \xi_i  \Big] \\
& \qquad = \frac{1}{b n} \left( \frac{n-b}{n-1} \right) \Big[ - n \Big\| \frac{1}{n}  \sum_{i=1}^{n} \xi_i \Big\|^2 +  \sum_{i=1}^{n} \| \xi_i \|^2  \Big] \\
& \qquad \leq \frac{1}{b} \left( \frac{n-b}{n-1} \right) \frac{1}{n}\sum_{i=1}^{n} \| \xi_i \|^2 \\
& \qquad \overset{\eqref{xi_value}}{=} \frac{1}{b} \left( \frac{n-b}{n-1} \right) \frac{1}{n}\sum_{i=1}^{n} \| \nabla f_{i} (w_{j}^{(s)}) - \nabla f_{i}(w_{j-1}^{(s)}) \|^2 \\
& \qquad \overset{\eqref{eq:Lsmooth_basic}}{\leq} \frac{1}{b} \left( \frac{n-b}{n-1} \right) L^2 \eta^2 \| v_{j-1}^{(s)} \|^2
\end{align*}
Hence, by taking expectation, we have
\begin{align*}
\mathbb{E}[\| v_j^{(s)} - v_{j-1}^{(s)} \|^2 ] - \mathbb{E}[ \| \nabla F(w_{j}^{(s)}) - \nabla F(w_{j-1}^{(s)}) \|^2] \leq \frac{1}{b} \left( \frac{n-b}{n-1} \right) L^2 \eta^2 \mathbb{E}[ \| v_{j-1}^{(s)} \|^2]. 
\end{align*}

By Lemma \ref{lem:var_diff_01}, for $t \geq 1$, 
\begin{align*}
\mathbb{E}[ \| \nabla F(w_{t}^{(s)}) - v_{t}^{(s)} \|^2 ] 
& = \sum_{j = 1}^{t} \mathbb{E}[ \| v_{j}^{(s)} - v_{j-1}^{(s)} \|^2 ]  
 - \sum_{j = 1}^{t} \mathbb{E}[ \| \nabla F(w_{j}^{(s)}) - \nabla F(w_{j-1}^{(s)}) \|^2 ] \\
 & \leq \frac{1}{b} \left( \frac{n-b}{n-1} \right) L^2 \eta^2 \sum_{j = 1}^{t} \mathbb{E}[ \| v_{j-1}^{(s)} \|^2].  
\end{align*}

Note that $\| \nabla F(w_{0}^{(s)}) - v_{0}^{(s)} \|^2 = 0$. Hence, by summing over $t = 0,\dots,m$ ($m \geq 1$), we have
\begin{align*}
\sum_{t=0}^{m} \mathbb{E}\| \nabla F(w_{t}^{(s)}) - v_{t}^{(s)} \|^2 \leq \frac{1}{b} \left( \frac{n-b}{n-1} \right) L^2 \eta^2 \Big[ m  \mathbb{E}\|v_{0}^{(s)} \|^2 + (m-1) \mathbb{E}\|v_{1}^{(s)} \|^2 + \dots + \mathbb{E}\|v_{m-1}^{(s)}\|^2 \Big ]. 
\end{align*}

By choosing $\eta \leq \frac{2}{L\left(\sqrt{1 + \frac{4m}{b}\left(\frac{n-b}{n-1} \right)} + 1\right)}$, we have
\begin{align*}
& \sum_{t=0}^{m} \mathbb{E}[ \| \nabla F(w_{t}^{(s)}) - v_{t}^{(s)} \|^2 ]  - ( 1 - L\eta ) \sum_{t=0}^{m} \mathbb{E} [ \| v_{t}^{(s)} \|^2 ] \\
& \leq \frac{1}{b} \left( \frac{n-b}{n-1} \right) L^2 \eta^2 \Big[ m  \mathbb{E}\|v_{0}^{(s)} \|^2 + (m-1) \mathbb{E}\|v_{1}^{(s)} \|^2 + \dots + \mathbb{E}\|v_{m-1}^{(s)}\|^2 \Big ] \\ &- (1 - L\eta) \Big [ \mathbb{E}\|v_{0}^{(s)} \|^2 + \mathbb{E}\|v_{1}^{(s)} \|^2 + \dots + \mathbb{E}\|v_{m}^{(s)}\|^2  \Big ] \\
& \leq \Big[\frac{1}{b} \left( \frac{n-b}{n-1} \right) L^2\eta^2 m - (1 - L\eta) \Big] \sum_{t=1}^{m} \mathbb{E} [ \| v_{t-1}^{(s)} \|^2 ] \leq 0,  \tagthis \label{eq:equal_zero_mb}
\end{align*}

since $\eta = \frac{2}{L\left(\sqrt{1 + \frac{4m}{b}\left(\frac{n-b}{n-1} \right)} + 1\right)}$ is a root of equation $\frac{1}{b} \left( \frac{n-b}{n-1} \right) L^2\eta^2 m - (1 - L\eta) = 0$. 

Therefore, with $\eta \leq \frac{2}{L\left(\sqrt{1 + \frac{4m}{b}\left(\frac{n-b}{n-1} \right)} + 1\right)}$, we have
\begin{align*}
\mathbb{E}[ F(w_{m+1}^{(s)}) ]  &\leq \mathbb{E} [ F(w_0^{(s)}) ] - \frac{\eta}{2} \sum_{t=0}^{m} \mathbb{E}[ \| \nabla F(w_{t}^{(s)})\|^2 ]. 
\end{align*}

Following the same derivation of Theorem~\ref{thm_main_nonconvex}, we could achieve the desired result as follows for any given $\tilde{w}_0$.  
\begin{align*}
    \frac{1}{(m+1)S}\sum_{s=1}^S \sum_{t=0}^{m} \mathbb{E}[ \| \nabla F(w_{t}^{(s)})\|^2 ] \leq \frac{2}{\eta [(m+1) S]} [ F(\tilde{w}_0) - F^*],
\end{align*}
where $F^*$ is any lower bound of $F$, and $w_{t}^{(s)}$ is the $t$-th iteration in the $s$-th outer loop. 
\end{proof}

\subsection*{Proof of Corollary \ref{cor_main_derivation_nonconvex_mb}}

\textbf{Corollary \ref{cor_main_derivation_nonconvex_mb}. }
\textit{For the conditions in Theorem~\ref{thm_main_nonconvex_mb}, in order to achieve an $\epsilon$-accurate solution, the total complexity is $$\Ocal\left(  \left[ \left( \frac{n + 2 b m}{m + 1}    \right)  \left( \sqrt{1 + \frac{4m}{b}\left(\frac{n-b}{n-1} \right)} \right) \frac{1}{\epsilon} \right] \vee \left[ n + 2 b m \right] \right).$$ }

\begin{proof}
By Theorem~\ref{thm_main_nonconvex_mb}, let
\begin{align*}
	\eta = \frac{2}{L\left(\sqrt{1 + \frac{4m}{b}\left(\frac{n-b}{n-1} \right)} + 1\right)}.
\end{align*}
Hence, we have
\begin{align*}
    \frac{1}{(m+1)S}\sum_{s=1}^S \sum_{t=0}^{m} \mathbb{E}[ \| \nabla F(w_{t}^{(s)})\|^2 ] & \leq \frac{2}{\eta [(m+1) S]} [ F(\tilde{w}_0) - F^*] \\
    	&= \frac{\left(\sqrt{1 + \frac{4m}{b}\left(\frac{n-b}{n-1} \right)} + 1\right)}{(m+1) S} L [ F(\tilde{w}_0) - F^*] \\
    	& \leq \frac{\left(\sqrt{1 + \frac{4m}{b}\left(\frac{n-b}{n-1} \right)}\right)}{(m+1) S} 2 L [ F(\tilde{w}_0) - F^*] = \epsilon.
\end{align*}
In order to achieve the $\epsilon$-accurate solution, we need
\begin{align*}
    S &= \frac{\left(\sqrt{1 + \frac{4m}{b}\left(\frac{n-b}{n-1} \right)}\right)}{(m+1) \epsilon} 2 L [ F(\tilde{w}_0) - F^*] = \Ocal\left( \frac{\left(\sqrt{1 + \frac{4m}{b}\left(\frac{n-b}{n-1} \right)}\right)}{(m+1) \epsilon} \vee 1  \right), 
\end{align*}
since $S \geq 1$. Therefore, the total complexity is
\begin{align*}
	(n + b \cdot 2 m) S & =  \Ocal\left(  \left[ \left( \frac{n + 2 b m}{m + 1}    \right)  \left( \sqrt{1 + \frac{4m}{b}\left(\frac{n-b}{n-1} \right)} \right) \frac{1}{\epsilon} \right] \vee \left[ n + 2 b m \right] \right). 
\end{align*}
\end{proof}

\subsection*{Proof of Corollary \ref{cor_main_derivation_nonconvex_2_mb}}

\textbf{Corollary \ref{cor_main_derivation_nonconvex_2_mb}. }
\textit{For the conditions in Theorem~\ref{thm_main_nonconvex_mb} and Corollary~\ref{cor_main_derivation_nonconvex_mb} with $b = n^{\alpha}$ and $m = n^{\beta}$ where $\alpha + \beta = 1$ with $\beta \geq 1/2$ and $0 \leq \alpha \leq 1/2$, in order to achieve an $\epsilon$-accurate solution, the total complexity is $$\Ocal\left( \frac{\sqrt{n}}{\epsilon} \vee n \right).$$ }
\begin{proof}
Let $b = n^{\alpha}$, $\alpha < 1$, and $m = n^{\beta}$, we have
\begin{align*}
    \left( \frac{n + 2 b m}{m + 1}    \right)  \sqrt{1 + \frac{4m}{b}\left(\frac{n-b}{n-1} \right)} &= \left( \frac{n + 2 n^{\alpha + \beta}}{n^{\beta} + 1}    \right)  \sqrt{1 + 4 n^{\beta - \alpha}\left(\frac{n- n^{\alpha}}{n-1} \right)} \\ & \leq \frac{n + 2 n^{\alpha + \beta} }{n^{\beta}} 2 \sqrt{1 + n^{\beta - \alpha}}.
\end{align*}
If $\beta \geq \alpha$, we have
\begin{align*}
    \frac{n + 2 n^{\alpha + \beta} }{n^{\beta}} 2 \sqrt{1 + n^{\beta - \alpha}} \leq 2\sqrt{2} \left( \frac{n + 2 n^{\alpha + \beta} }{n^{\beta}} \right) n^{(\beta - \alpha)/2} = 2\sqrt{2} ( n^{1 - \alpha/2 - \beta/2} + 2 n^{\alpha/2 + \beta/2} ). 
\end{align*}
In order to minimize the order of $n$, we need to choose $1 - \alpha/2 - \beta/2 = \alpha/2 + \beta/2$, which is equivalent to $\alpha + \beta = 1$ with $\beta \geq \alpha$. The best option is to choose $\alpha + \beta = 1$ with $\beta \geq 1/2$ and $0 \leq \alpha \leq 1/2$ in order to achieve $\Ocal(n^{1/2})$. 

If $\beta \leq \alpha$, we have
\begin{align*}
    \frac{n + 2 n^{\alpha + \beta} }{n^{\beta}} 2 \sqrt{1 + n^{\beta - \alpha}} \leq 2\sqrt{2} (n^{1-\beta} + 2 n^{\alpha} ).
\end{align*}
In order to minimize the order of $n$, we need to choose $1 - \beta = \alpha$, which is equivalent to $\alpha + \beta = 1$ with $\beta \leq \alpha$. The best option is to choose $\beta = 1/2$ and $\alpha = 1/2$ in order to achieve $\Ocal(n^{1/2})$. 

Therefore, with $b = n^{\alpha}$ and $m = n^{\beta}$ where $\alpha + \beta = 1$ with $\beta \geq 1/2$ and $0 \leq \alpha \leq 1/2$, we have
\begin{align*}
    \left( \frac{n + 2 b m}{m + 1}    \right)  \sqrt{1 + \frac{4m}{b}\left(\frac{n-b}{n-1} \right)} & = \Ocal(n^{1/2}). 
\end{align*}
By Corollary~\ref{cor_main_derivation_nonconvex_mb} with $b m = n^{\alpha + \beta} = n$, it implies the total complexity
\begin{align*}
    (n + b \cdot 2 m) S & =  \Ocal\left( \frac{\sqrt{n}}{\epsilon} \vee n \right). 
\end{align*}
\end{proof}

\section*{Convex SARAH++}

\subsection*{Proof of Lemma \ref{lem_basic_lem_01}}

\textbf{Lemma \ref{lem_basic_lem_01}}. \textit{Suppose that Assumptions \ref{ass_Lsmooth} and \ref{ass_convex} holds. Consider SARAH (Algorithm~\ref{sarah}) within a single outer loop with $\eta \leq \frac{1}{L}$. Then, for $t \geq 0$ and any $s \geq 1$, we have
\begin{align*}
    \mathbb{E}[F(w_{t+1}^{(s)}) - F(w_*)] \leq \mathbb{E}[F(w_{t}^{(s)}) - F(w_*)]  - \frac{\eta}{2} \mathbb{E} [ \| \nabla F(w_t^{(s)}) \|^2]  + \frac{\eta}{2} \left( L\eta \mathbb{E}[\| v_0^{(s)} \|^2] - \mathbb{E}[\| v_t^{(s)} \|^2] \right), 
\end{align*}
where $w_*$ is any optimal solution of $F$.}

\begin{proof}
By using \eqref{eq_proof_lem_000} and adding $-F(w_*)$ for both sides, where $w_* = \arg \min_{w} F(w)$, we have
\allowdisplaybreaks
\begin{align*}
\mathbb{E}[ F(w_{t+1}^{(s)}) - F(w_*)] & \leq \mathbb{E}[ F(w_{t}^{(s)}) - F(w_*)] - \frac{\eta}{2} \mathbb{E}[ \| \nabla F(w_{t}^{(s)})\|^2 ] 
+ \frac{\eta}{2} \mathbb{E}[ \| \nabla F(w_{t}^{(s)}) - v_{t}^{(s)} \|^2 ] \\ & \qquad - \left( \frac{\eta}{2} - \frac{L\eta^2}{2} \right) \mathbb{E} [ \| v_{t}^{(s)} \|^2 ] \\
& \overset{\eqref{eq:bound1}}{\leq} \mathbb{E}[ F(w_{t}^{(s)}) - F(w_*)] - \frac{\eta}{2} \mathbb{E}[ \| \nabla F(w_{t}^{(s)})\|^2 ] \\
& \qquad + \frac{\eta}{2} \frac{\eta L}{(2 - \eta L)} \left( \mathbb{E}[ \|v_{0}^{(s)} \|^2] - \mathbb{E}[\| v_{t}^{(s)} \|^2] \right) - \left( \frac{\eta}{2} - \frac{L\eta^2}{2} \right) \mathbb{E} [ \| v_{t}^{(s)} \|^2 ] \\
& = \mathbb{E}[ F(w_{t}^{(s)}) - F(w_*)] - \frac{\eta}{2} \mathbb{E}[ \| \nabla F(w_{t}^{(s)})\|^2 ] \\
& \qquad + \frac{\eta}{2} \left( \frac{\eta L}{(2 - \eta L)} \left( \mathbb{E}[ \|v_{0}^{(s)} \|^2] - \mathbb{E}[\| v_{t}^{(s)} \|^2] \right) - \left( 1 - L\eta \right) \mathbb{E} [ \| v_{t}^{(s)} \|^2 ] \right) \\
& \overset{\eta\leq\frac{1}{L}}{\leq} \mathbb{E}[ F(w_{t}^{(s)}) - F(w_*)] - \frac{\eta}{2} \mathbb{E}[ \| \nabla F(w_{t}^{(s)})\|^2 ] \\
& \qquad + \frac{\eta}{2} \left( \eta L \left( \mathbb{E}[ \|v_{0}^{(s)} \|^2] - \mathbb{E}[\| v_{t}^{(s)} \|^2] \right) - \left( 1 - L\eta \right) \mathbb{E} [ \| v_{t}^{(s)} \|^2 ] \right) \\
& = \mathbb{E}[ F(w_{t}^{(s)}) - F(w_*)] - \frac{\eta}{2} \mathbb{E}[ \| \nabla F(w_{t}^{(s)})\|^2 ] + \frac{\eta}{2} \left( L\eta \mathbb{E}[\| v_0^{(s)} \|^2] - \mathbb{E}[\| v_t^{(s)} \|^2] \right). 
\end{align*}
\end{proof}

\subsection*{Proof of Theorem \ref{thm_general_convex_sarahpp}}

\textbf{Theorem~\ref{thm_general_convex_sarahpp}} (Smooth general convex). \textit{Suppose that Assumptions \ref{ass_Lsmooth} and \ref{ass_convex} holds. Consider SARAH++ (Algorithm~\ref{sarah_new}) with $\eta \leq \frac{\gamma}{L}$, $0 < \gamma \leq 1$. Then, the expectation of the average of squared norm of gradient of all iterations generated by SARAH++}
\begin{align*}
    \mathbb{E}\left[ \frac{1}{T_1 + \dots + T_S}  \sum_{s=1}^S \sum_{t=0}^{T_s-1} \mathbb{E} [ \| \nabla F(w_t^{(s)}) \|^2 |T_1,\dots,T_S ]  \right] \leq \frac{2}{T \eta} [ F(\tilde{w}_0) - F(w_*) ].
\end{align*}

\begin{proof}
We recall the following definitions. $T_s$ is the stopping time (a random variable) of the $s$-th outer iteration such that
\begin{align*}
    T_s = \min \left\{ \min_{t \geq 0} \left\{ t : \| v_{t}^{(s)} \|^2 < \gamma \| v_0^{(s)} \|^2 \right\} , m + 1 \right\} \ , \ s = 1,2,\dots
\end{align*}
and $S$ is the stopping time of the outer iterations (a random variable) and such that for some $T > 0$
\begin{align*}
    S = \min_{\hat{S}} \left\{ \hat{S} : \sum_{s=1}^{\hat{S}} T_s \geq T  \right\}. 
\end{align*}


Note that $T_s \geq 1$ is the first time such that $\| v_{T_s}^{(s)} \|^2 < \gamma \| v_0^{(s)} \|^2$. Hence, for a given $T_s$, we have $\| v_{t}^{(s)} \|^2 \geq \gamma \| v_0^{(s)} \|^2$, for $0 \leq t \leq T_s-1$, and 
\begin{align*}
    \mathbb{E}[F(w_{{T_s}}^{(s)}) - F(w_*)] & \leq \mathbb{E}[F(w_{{{T_s}-1}}^{(s)}) - F(w_*)]  - \frac{\eta}{2} \mathbb{E} [ \| \nabla F(w_{{T_s}-1}^{(s)}) \|^2] + \frac{\eta}{2} \left( L\eta \mathbb{E}[\| v_0^{(s)} \|^2] - \mathbb{E}[\| v_{{T_s}-1}^{(s)} \|^2] \right) \\
    & \overset{\eta \leq \frac{\gamma}{L}}{\leq} \mathbb{E}[F(w_{{{T_s}-1}}^{(s)}) - F(w_*)]  - \frac{\eta}{2} \mathbb{E} [ \| \nabla F(w_{{T_s}-1}^{(s)}) \|^2] + \frac{\eta}{2} \left( \gamma \mathbb{E}[\| v_0^{(s)} \|^2] - \mathbb{E}[\| v_{{T_s}-1}^{(s)} \|^2] \right) \\
    & \leq \mathbb{E}[F(w_{{{T_s}-1}}^{(s)}) - F(w_*)]  - \frac{\eta}{2} \mathbb{E} [ \| \nabla F(w_{{T_s}-1}^{(s)}) \|^2] \\
    & \leq \mathbb{E}[F(w_{0}^{(s)}) - F(w_*)] - \frac{\eta}{2} \sum_{t=0}^{T_s-1} \mathbb{E} [ \| \nabla F(w_t^{(s)}) \|^2  ]. 
\end{align*}
Since $\tilde{w}_{s} = w_{T_s}^{(s)}$ and $\tilde{w}_{s-1} = w_0^{(s)}$, for given $T_1,\dots,T_S$, we have
\begin{align*}
    \mathbb{E}[F(\tilde{w}_{S}) - F(w_*)] & \leq \mathbb{E}[F(\tilde{w}_{S-1}) - F(w_*)] - \frac{\eta}{2} \sum_{t=0}^{T_S-1} \mathbb{E} [ \| \nabla F(w_t^{(s)}) \|^2  ] \\
    & \leq \mathbb{E}[F(\tilde{w}_{0}) - F(w_*)] - \frac{\eta}{2} \sum_{s=1}^S \sum_{t=0}^{T_s-1} \mathbb{E} [ \| \nabla F(w_t^{(s)}) \|^2  ]. 
\end{align*}

Since $F(\tilde{w}_{S}) \geq F(w_*)$, bringing the second term of the RHS to the LHS. For any given $\tilde{w}_{0}$, we have
\begin{align*}
    \frac{\eta}{2} \sum_{s=1}^S \sum_{t=0}^{T_s-1} \mathbb{E} [ \| \nabla F(w_t^{(s)}) \|^2 |T_1,\dots,T_S ] \leq [F(\tilde{w}_{0}) - F(w_*)], 
\end{align*}
which is equivalent to
\begin{align*}
   \frac{1}{T_1 + \dots + T_S}  \sum_{s=1}^S \sum_{t=0}^{T_s-1} \mathbb{E} [ \| \nabla F(w_t^{(s)}) \|^2 |T_1,\dots,T_S ] & \leq \frac{1}{T_1 + \dots + T_S} \frac{2}{\eta} [F(\tilde{w}_{0}) - F(w_*)] \\
   & \leq \frac{2}{\eta T} [F(\tilde{w}_{0}) - F(w_*)],  
\end{align*}
where the last inequality follows since $\sum_{s=1}^S T_s \geq T$. Hence, by taking the expectation to both sides, we could have
\begin{align*}
    \mathbb{E}\left[ \frac{1}{T_1 + \dots + T_S}  \sum_{s=1}^S \sum_{t=0}^{T_s-1} \mathbb{E} [ \| \nabla F(w_t^{(s)}) \|^2 |T_1,\dots,T_S ]  \right] \leq \frac{2}{\eta T} [F(\tilde{w}_{0}) - F(w_*)]. 
\end{align*}
Therefore, we achieve the desired result since the LHS is the expectation of the average of squared norm of gradient of all iterations generated by SARAH++ (Algorithm~\ref{sarah_new}). 
\end{proof}

\subsection*{Proof of Corollary \ref{cor_general_convex_sarahpp}}

\textbf{Corollary~\ref{cor_general_convex_sarahpp}} (Smooth general convex). \textit{Consider the conditions in Theorem~\ref{thm_general_convex_sarahpp} with $\eta = \Ocal(\frac{1}{L})$. Then we could achieve the $\epsilon$-accurate solution after $\Ocal(\frac{1}{\epsilon})$ total iterations.}
\begin{proof}
The proof is trivial since we want
\begin{align*}
    \frac{2}{\eta T} [F(\tilde{w}_{0}) - F(w_*)] = \epsilon, 
\end{align*}
which requires $T = \frac{2[F(\tilde{w}_{0}) - F(w_*)]}{\eta} \cdot \frac{1}{\epsilon} = \Ocal(\frac{1}{\epsilon})$ iterations, where we could choose $\eta = \Ocal(\frac{1}{L})$. 
\end{proof}

\subsection*{Proof of Theorem \ref{thm_strongly_convex_sarahpp}}

\textbf{Theorem~\ref{thm_strongly_convex_sarahpp}} (Smooth strongly convex). \textit{Suppose that Assumptions \ref{ass_Lsmooth}, \ref{ass_stronglyconvex} and \ref{ass_convex} holds. Consider SARAH++ (Algorithm~\ref{sarah_new}) with $\eta \leq \frac{\gamma}{L}$, $0 < \gamma \leq 1$. Then, for the final output $\hat{w}$ of SARAH++, we have}
\begin{align*}
    \mathbb{E} [ F(\hat{w}) - F(w_*) ] &\leq (1 - \mu \eta)^T [ F(\tilde{w}_0) - F(w_*) ]. 
\end{align*}

\begin{proof}
Following the beginning part of the proof of Theorem~\ref{thm_general_convex_sarahpp}, we have, for a given $T_s$, 
\begin{align*}
    \mathbb{E}[F(w_{{T_s}}^{(s)}) - F(w_*)] & \leq \mathbb{E}[F(w_{{{T_s}-1}}^{(s)}) - F(w_*)]  - \frac{\eta}{2} \mathbb{E} [ \| \nabla F(w_{{T_s}-1}^{(s)}) \|^2] \\
    & \overset{\eqref{eq:strongconvexity2}}{\leq} (1 - \mu \eta) \mathbb{E}[F(w_{{{T_s}-1}}^{(s)}) - F(w_*)] \\
    & \leq (1 - \mu \eta)^{T_s} \mathbb{E}[F(w_{0}^{(s)}) - F(w_*)]
\end{align*}
Since $\tilde{w}_{s} = w_{T_s}^{(s)}$ and $\tilde{w}_{s-1} = w_0^{(s)}$, for given $T_1,\dots,T_S$, we have
\begin{align*}
    \mathbb{E} [ F(\hat{w}) - F(w_*) | T_1,\dots,T_S ] & = \mathbb{E} [ F(\tilde{w}_S) - F(w_*) | T_1,\dots,T_S ] \\ 
    & \leq (1 - \mu \eta)^{T_1 + \dots + T_S} [ F(\tilde{w}_0) - F(w_*) ] \\
    & \leq (1 - \mu \eta)^T [ F(\tilde{w}_0) - F(w_*) ], 
\end{align*}
where the last inequality follows since $\sum_{s=1}^S T_s \geq T$. Hence, by taking the expectation to both sides, we could have
\begin{align*}
    \mathbb{E} [ F(\hat{w}) - F(w_*)] &\leq (1 - \mu \eta)^T [ F(\tilde{w}_0) - F(w_*) ]. 
\end{align*}
\end{proof}

\subsection*{Proof of Corollary \ref{cor_strongly_convex_sarahpp}}

\textbf{Corollary~\ref{cor_strongly_convex_sarahpp}} (Smooth strongly convex). \textit{Consider the conditions in Theorem~\ref{thm_strongly_convex_sarahpp} with $\eta = \Ocal(\frac{1}{L})$. Then we could achieve $\mathbb{E} [ F(\hat{w}) - F(w_*) ] \leq \epsilon$ after $\Ocal(\kappa \log(\frac{1}{\epsilon}))$ total iterations, where $\kappa = L/\mu$ is the condition number.}

\begin{proof}
We want
\begin{align*}
    (1 - \mu \eta)^T [ F(\tilde{w}_0) - F(w_*) ] = \epsilon. 
\end{align*}
Hence, 
\begin{align*}
T = -\frac{1}{\log(1 - \mu\eta)} \log\left( \frac{[ F(\tilde{w}_{0}) - F(w_*) ]}{\epsilon} \right). 
\end{align*}
Note that: $-\frac{1}{x} - 1 \leq - \frac{1}{\log(1+x)} \leq -\frac{1}{x}$, $-1 < x < 0$. We can have
\begin{align*}
\left(\frac{1}{\mu \eta} - 1 \right) \log\left( \frac{[ F(\tilde{w}_{0}) - F(w_*) ]}{\epsilon} \right) \leq T \leq \frac{1}{\mu \eta} \log\left( \frac{[ F(\tilde{w}_{0}) - F(w_*) ]}{\epsilon} \right). 
\end{align*}
By choosing $\eta = \Ocal(\frac{1}{L})$, we have $T = \Ocal(\kappa \log(\frac{1}{\epsilon}))$. 
\end{proof}




\section*{Additional Experiments}\label{sec_additional_experiment}

We provide more experiments in this section on popular data sets with diverse size $n$ including \textit{covtype} ($n = 406,708$ training data; estimated $L \simeq 1.90$), \textit{ijcnn1} ($n = 91,701$ training data; estimated $L \simeq 1.77$), \textit{w8a} ($n = 49,749$ training data, estimated $L \simeq 7.05$) and \textit{phishing} ($n = 7,738$ training data, estimated  $L \simeq 7.49$) from LIBSVM. 



\subsection*{Additional experiments in Section \ref{sec_experiments}}

\begin{figure}[h]
 \centering
 \includegraphics[width=0.245\textwidth]{Figs/covtype_sarahpp_sarah_03.png}
 \includegraphics[width=0.245\textwidth]{Figs/ijcnn1_sarahpp_sarah_03.png}
 \includegraphics[width=0.245\textwidth]{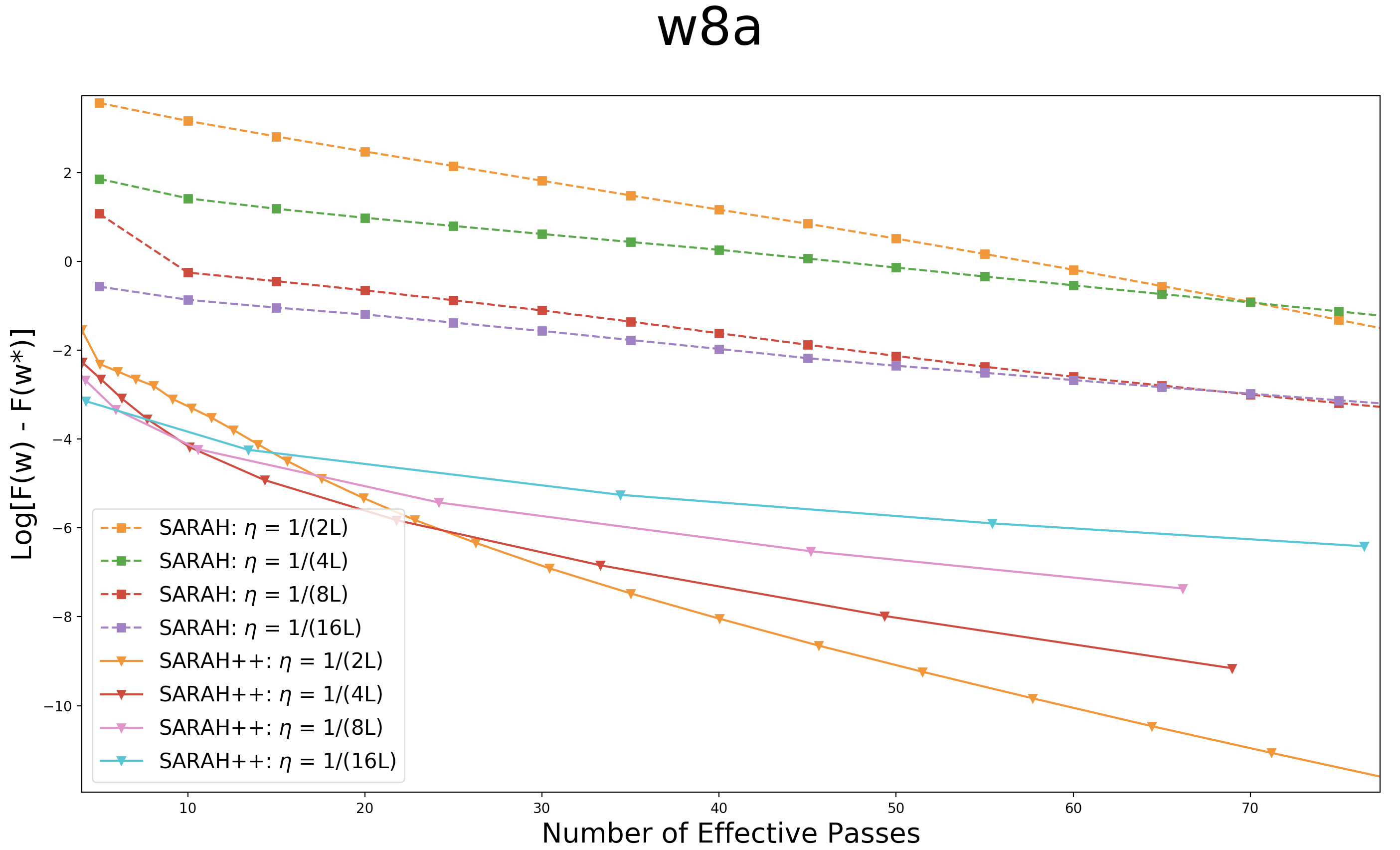} 
 \includegraphics[width=0.245\textwidth]{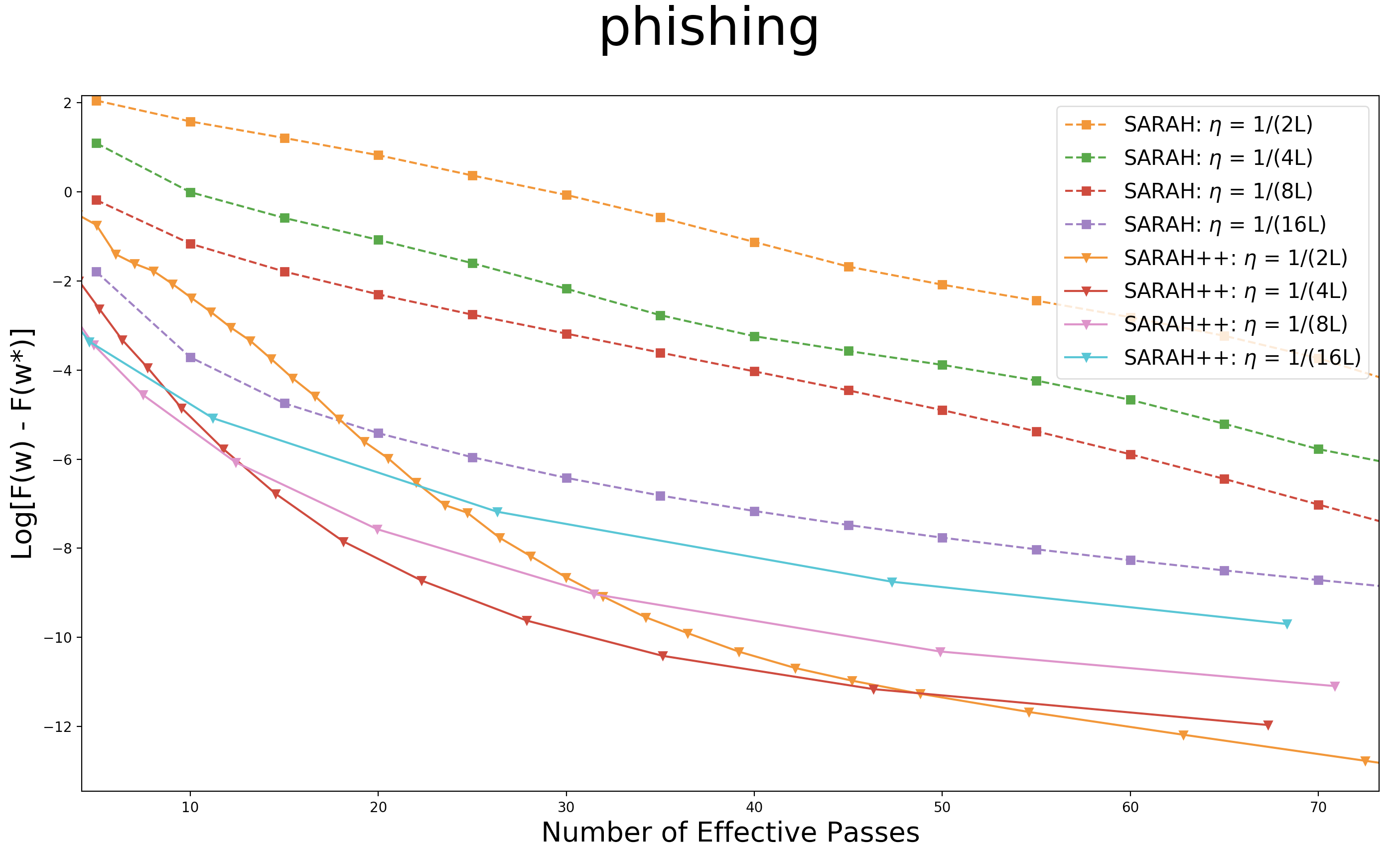} 
  \caption{Comparisons of $\log[F(w) - F(w_*)]$ between SARAH++ and SARAH with different learning rates on \textit{covtype}, \textit{ijcnn1}, \textit{w8a}, and \textit{phishing} datasets}
  \label{fig_sarahpp_sarah_apd}
 \end{figure}
 
 \subsection*{Additional experiments in Section \ref{sec_adaptive}}
 
\begin{figure}[h]
 \centering
 \includegraphics[width=0.245\textwidth]{Figs/covtype_sarah_adap_04.png}
 \includegraphics[width=0.245\textwidth]{Figs/ijcnn1_sarah_adap_04.png}
 \includegraphics[width=0.245\textwidth]{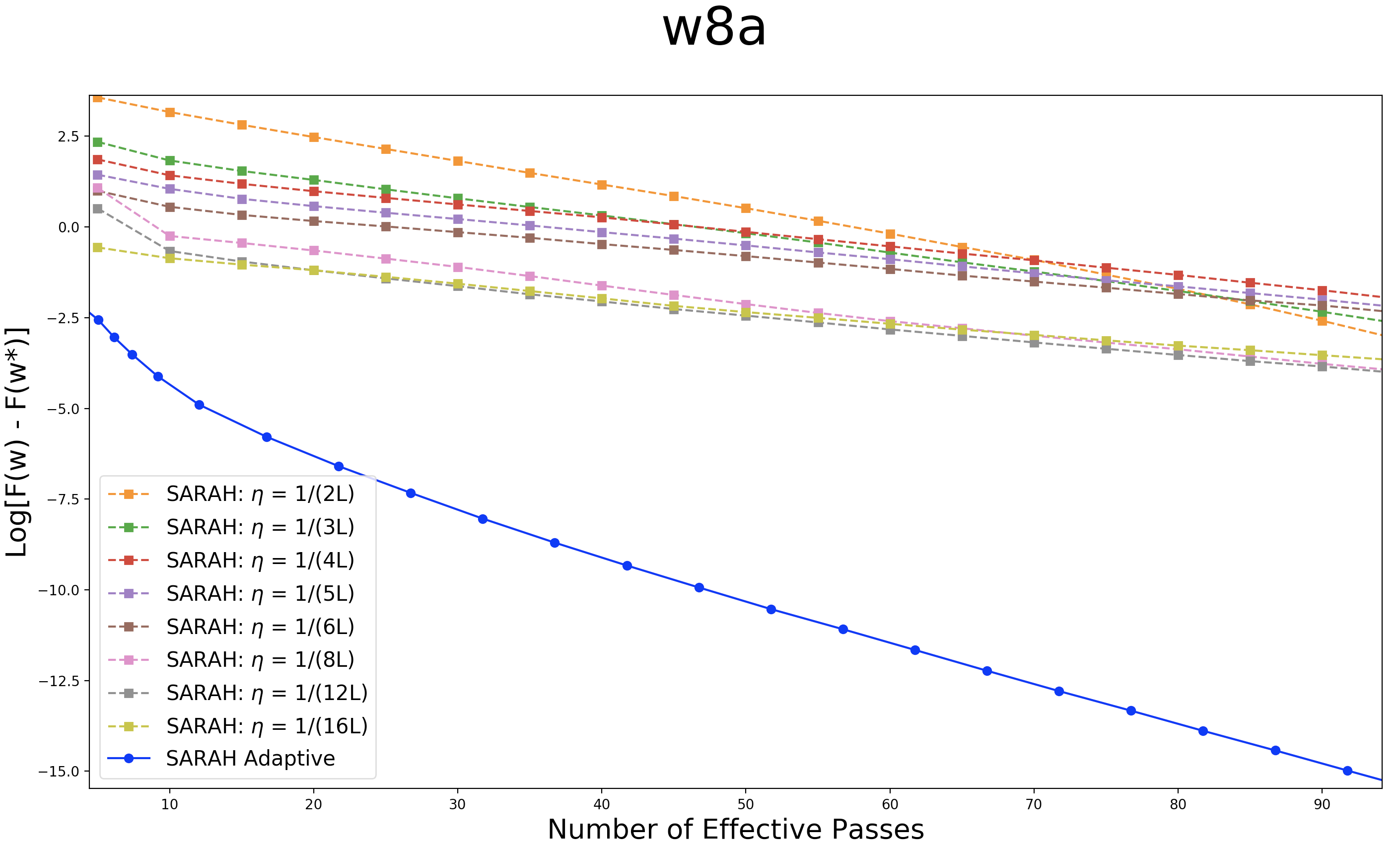} 
 \includegraphics[width=0.245\textwidth]{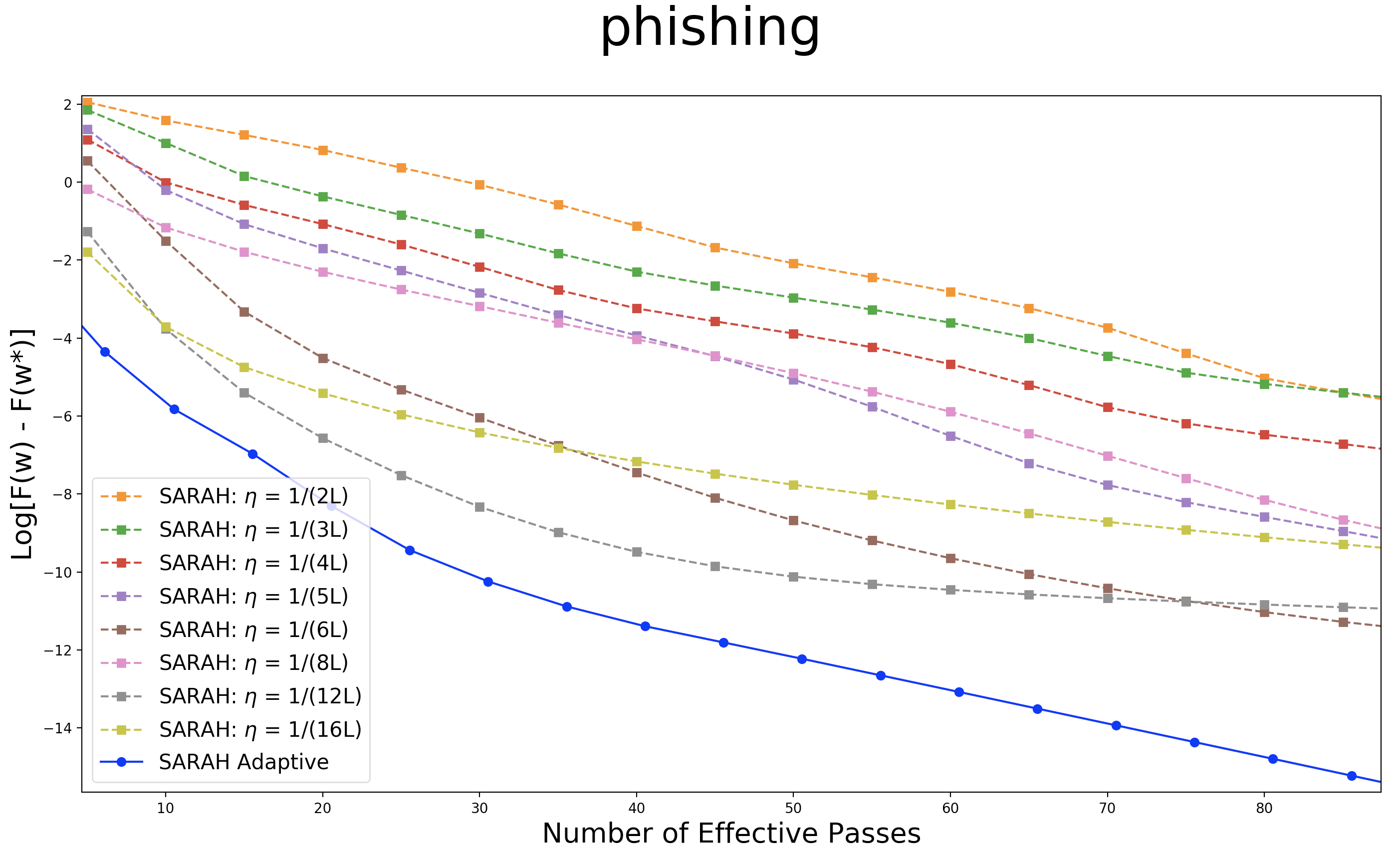} 
  \caption{Comparisons of $\log[F(w) - F(w_*)]$ between SARAH Adaptive and SARAH with different learning rates on \textit{covtype}, \textit{ijcnn1}, \textit{w8a}, and \textit{phishing} datasets}
  \label{fig_sarah_adap_sarah_apd}
 \end{figure}

\begin{figure}[h]
 \centering
 \includegraphics[width=0.245\textwidth]{Figs/covtype_sarah_adap_03.png}
 \includegraphics[width=0.245\textwidth]{Figs/ijcnn1_sarah_adap_03.png}
 \includegraphics[width=0.245\textwidth]{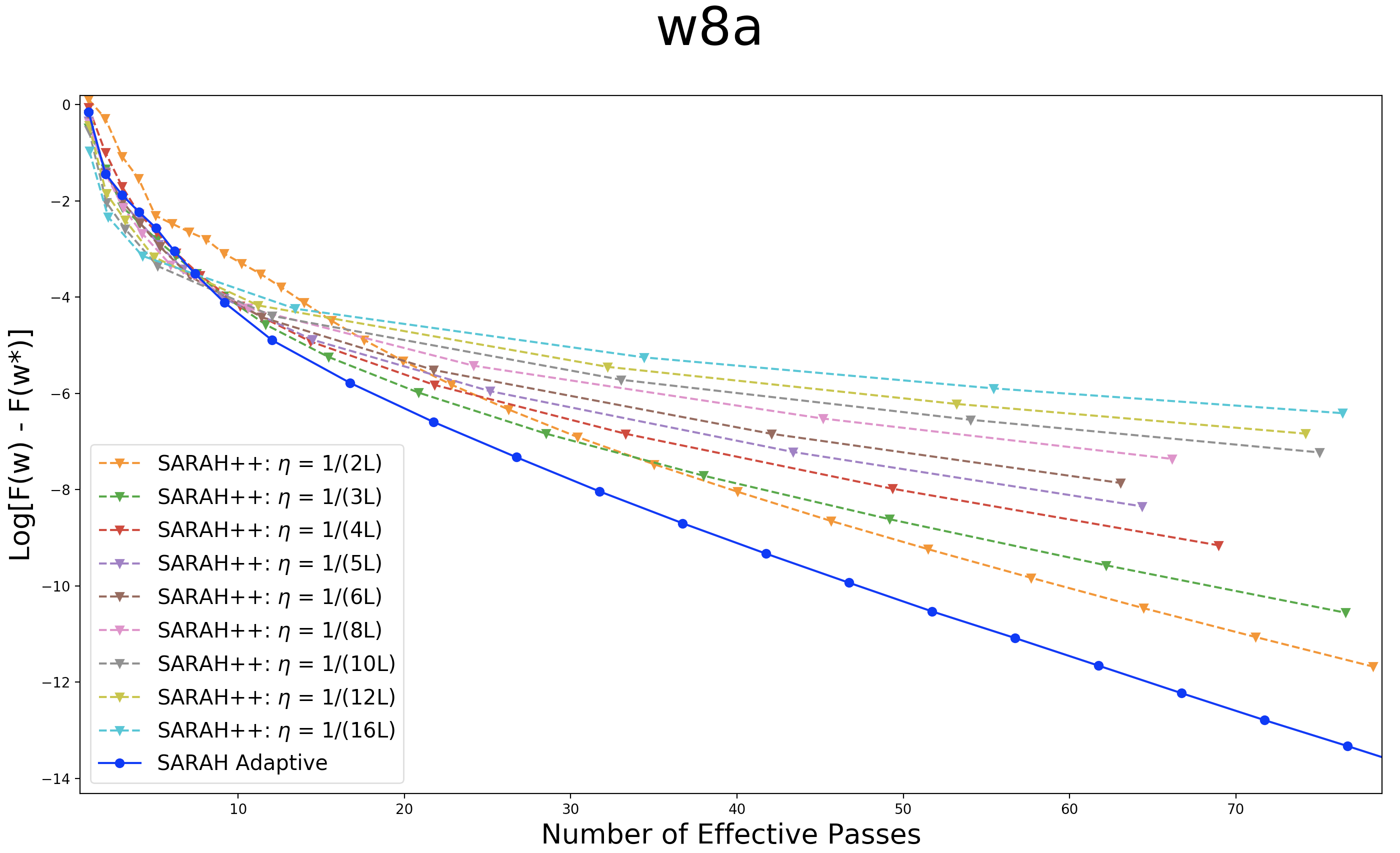} 
 \includegraphics[width=0.245\textwidth]{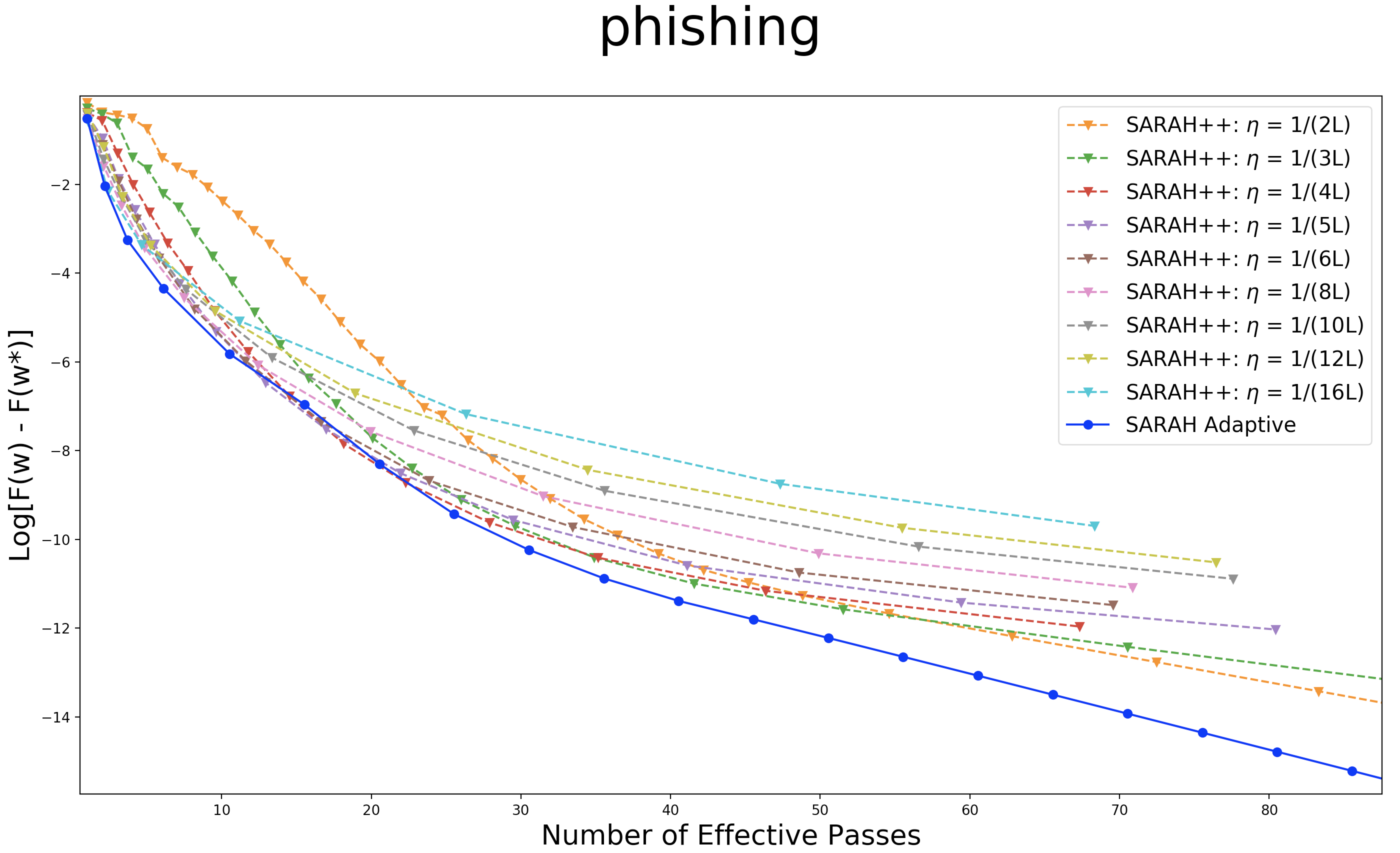} 
  \caption{Comparisons of $\log[F(w) - F(w_*)]$ between SARAH Adaptive and SARAH++ with different learning rates on \textit{covtype}, \textit{ijcnn1}, \textit{w8a}, and \textit{phishing} data sets}
  \label{fig_sarah_adap_sarahpp_apd}
 \end{figure}
 
\subsection*{Sensitivity of $\gamma$ for SARAH Adaptive}
 
\begin{figure}[h!]
 \centering
 \includegraphics[width=0.245\textwidth]{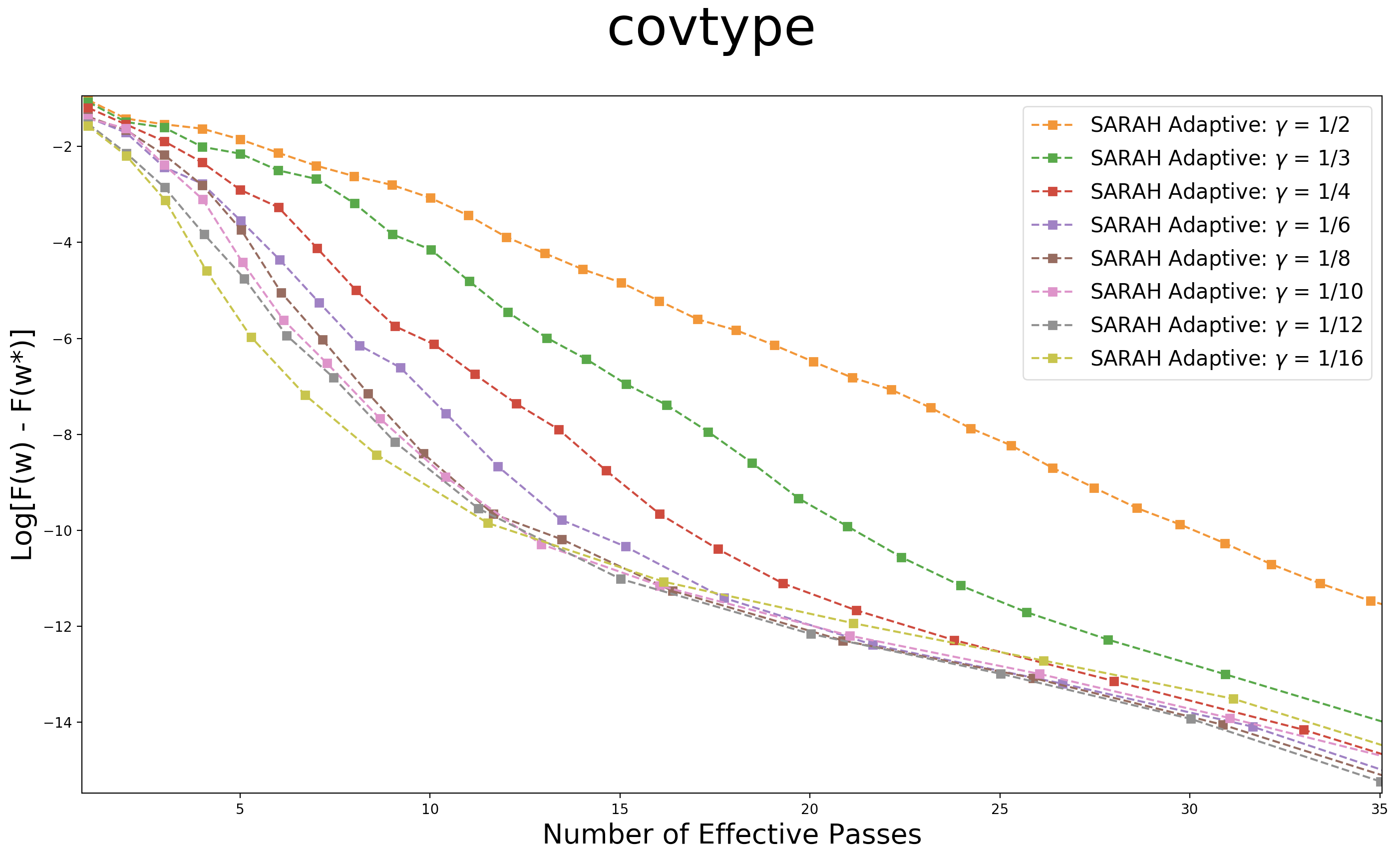}
 \includegraphics[width=0.245\textwidth]{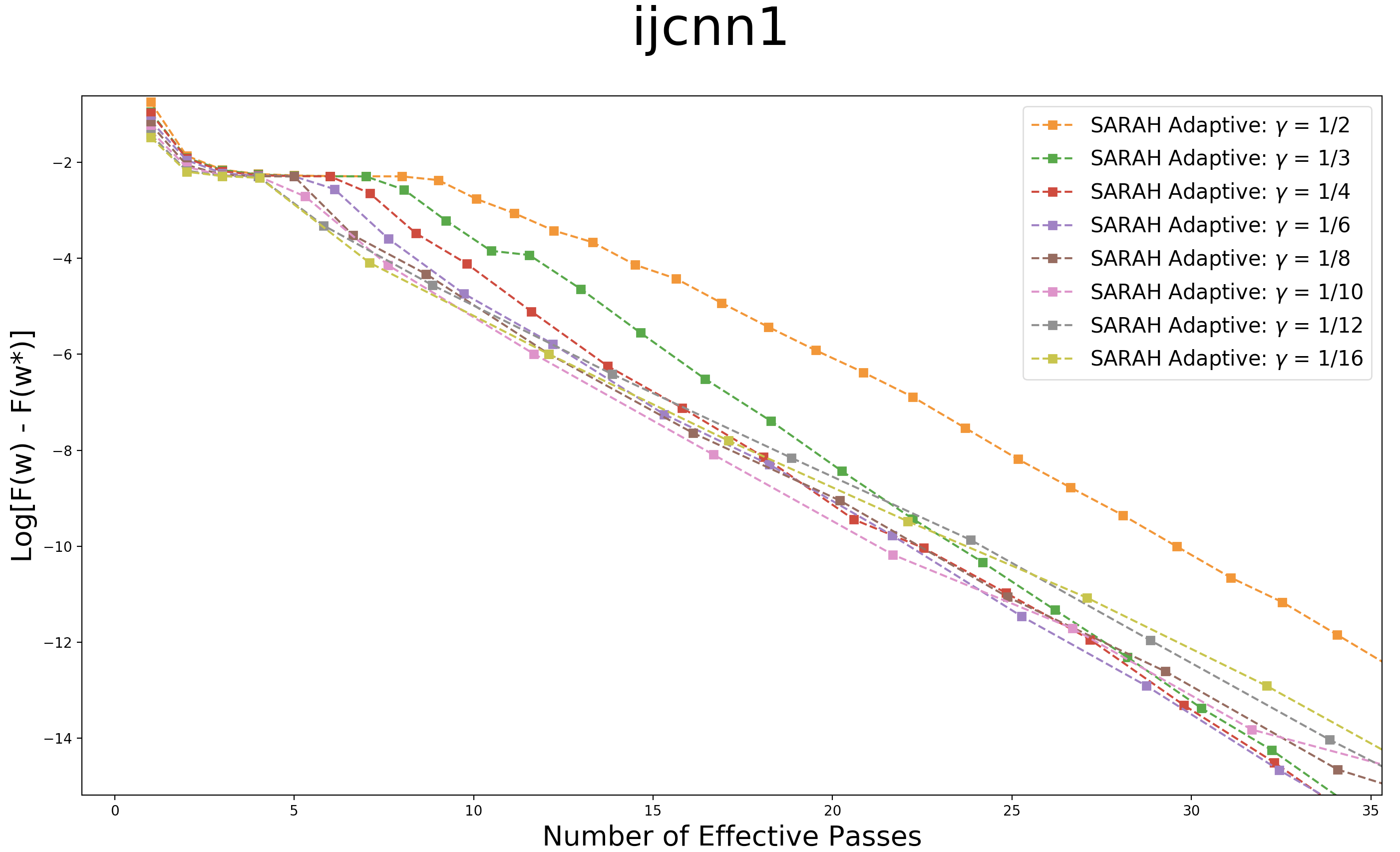}
 \includegraphics[width=0.245\textwidth]{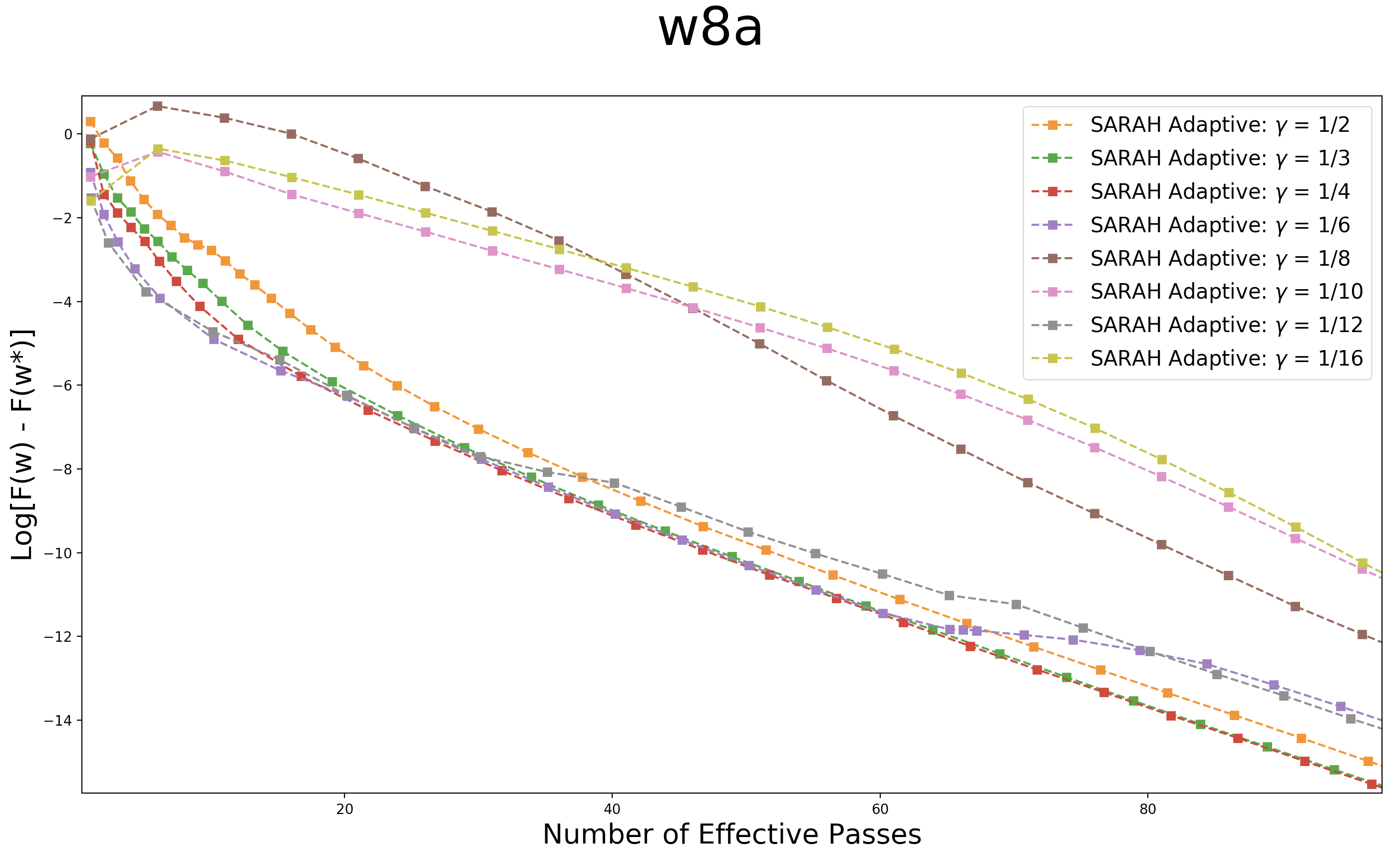} 
 \includegraphics[width=0.245\textwidth]{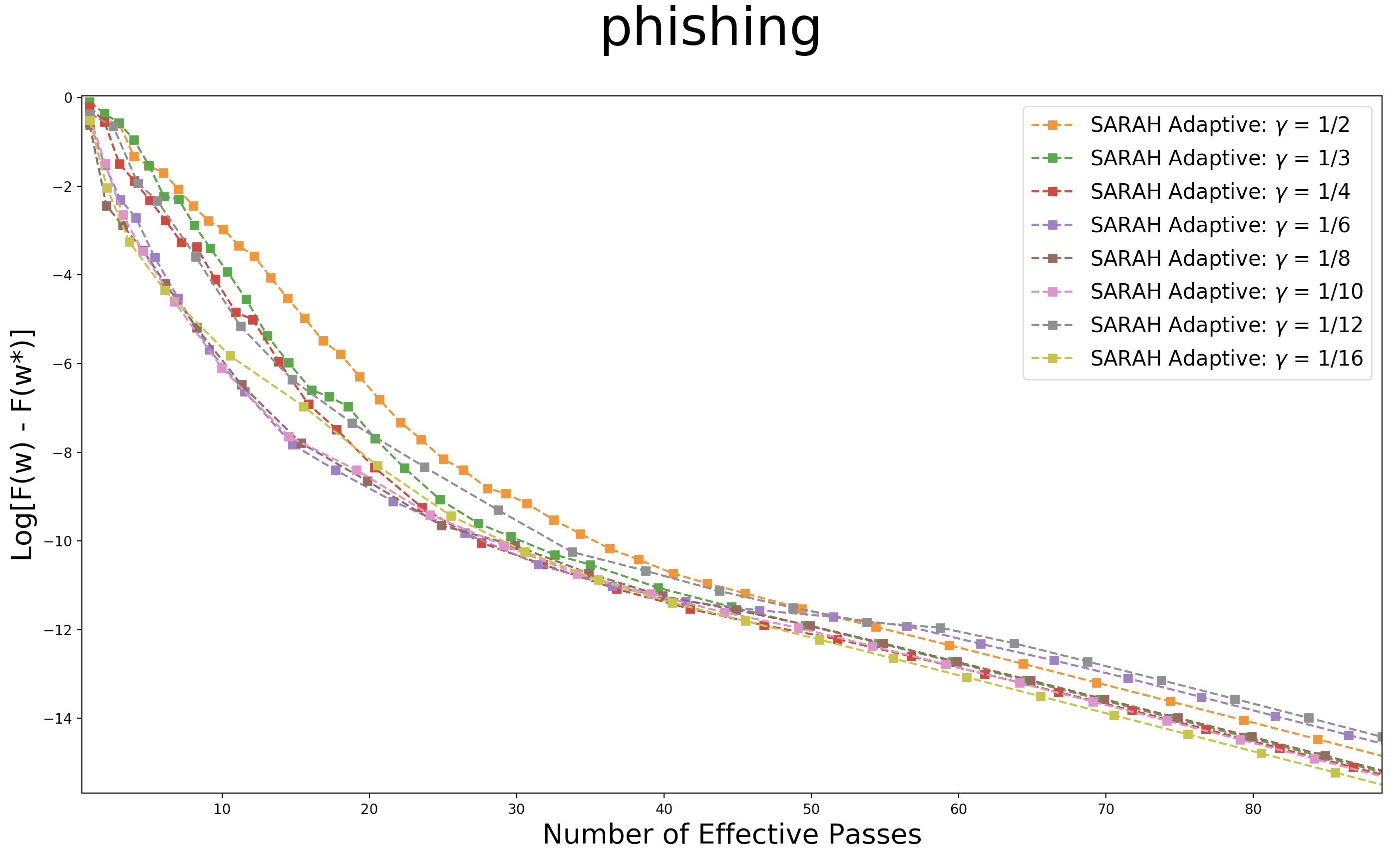} 
  \caption{Comparisons of $\log[F(w) - F(w_*)]$ with different value of $\gamma$ for SARAH Adaptive on \textit{covtype}, \textit{ijcnn1}, \textit{w8a}, and \textit{phishing} datasets}
  \label{fig_sarah_adap_gamma_apd}
 \end{figure}

In Figure~\ref{fig_sarah_adap_gamma_apd} we present the numerical performance of SARAH Adaptive for different values of \\ $\gamma = \left\{ \frac{1}{2}, \frac{1}{3}, \frac{1}{4}, \frac{1}{6}, \frac{1}{8}, \frac{1}{10}, \frac{1}{12}, \frac{1}{16}  \right\}$ on \textit{covtype}, \textit{ijcnn1}, \textit{w8a}, and \textit{phishing} data sets.  

 \clearpage